\documentclass[3p,times]{elsarticle}
\graphicspath{{Figures/}}
\DeclareGraphicsExtensions{.pdf,.png,.jpg}
\usepackage{amssymb}
\usepackage{amsthm}
\usepackage{graphicx}
\usepackage{amsmath,amsthm,amsbsy}
\usepackage{calligra,calrsfs,mathrsfs}
\usepackage{multirow}
\usepackage{cases}
\usepackage{pgfplots}
\pgfplotsset{compat=1.16}
\usepackage{subcaption}
\usepackage{tikz}
\usepackage{tkz-euclide}
\usetikzlibrary[patterns]
\usepackage{hyperref}
\usepackage{dsfont}
\usepackage{caption}
\usepackage{epstopdf}
\usepackage{epsfig}
\usepackage{array}
\newcolumntype{P}[1]{>{\centering\arraybackslash}p{#1}}
\usepackage{enumerate}
\usepackage{float}
\usepackage[T1]{fontenc}
\usepackage[latin9]{inputenc}
\usepgfplotslibrary{fillbetween}
\usetikzlibrary{patterns}
\usepackage{tikz-dimline}
\usepackage{color}
\usepackage{stmaryrd}
\usepackage{setspace}
\usepackage{amsthm}
\usepackage{textcomp}
\usepackage{soul,xcolor}
\usepackage{multirow}
\usepackage[inline]{enumitem}

\biboptions{sort&compress}
\usepackage{dirtytalk}
\usepackage{csvsimple}

\pgfplotsset{
/pgfplots/bar cycle list/.style={/pgfplots/cycle list={%
    {blue,fill=blue!30!white,mark=none},%
    {red,fill=red!30!white,mark=none},%
    {brown!60!black,fill=brown!30!white,mark=none},%
    {black,fill=gray,mark=none},%
    {violet!80!black,fill=violet,mark=none},%
    {orange,fill=orange!50!white,mark=none}%
    }
},
}
\usepackage[noend]{algpseudocode}
\usepackage{algorithm}

\usetikzlibrary{shapes.geometric} 
\usepackage[export]{adjustbox}

\usepackage{hyperref}
\hypersetup{
    colorlinks,
    linkcolor={magenta},
    citecolor={blue},
    urlcolor={blue!80!black},
    breaklinks=true,
	plainpages=true
}
    
   \usetikzlibrary{external}
\newcommand{\ifcomment}{\iffalse}
\newcommand{\bs}[1]{\boldsymbol{#1}}

\newcommand{\avg  }[1]{\langle #1 \rangle}

\newcommand{\ti}[1]{\tilde{#1}}
\newcommand{\G}{\Gamma}

\newcommand{\Om}{\Omega}

\newcommand{\cT}{{\mathcal T}}
\newcommand{\tO}{\tilde{\Omega}_h}
\newcommand{\hO}{\hat{\Omega}_h}
\newcommand{\Oe}{\Omega_{\mathrm{ext}}}

\newcommand{\tG}{\tilde{\Gamma}_h}

\newcommand{\tGD}{\tilde{\Gamma}_{D,h}}

\newcommand{\tT}{{\tilde{\cT}_h}}
\newcommand{\tS}{\mathrm{S}_{h}}
\newcommand{\tSLE}{\mathbf{S}_{h}}

\newtheorem{thm}{Theorem}

\newtheorem{prop}{Proposition}

\newdefinition{rem}{Remark}

\newtheorem*{remark}{Remark}

\newtheorem{definition}{Definition}[section]

\newcommand{\cref}[2]{\hyperref[#2]{#1~\ref*{#2}}}
\newcommand{\colref}[2]{\hyperref[#2]{#1~\ref*{#2}}}
\newcommand{\eqnref}[1]{\colref{Eq.}{#1}}
\newcommand{\figref}[1]{\colref{Figure}{#1}}

\newcommand{\secref}[1]{\colref{Section}{#1}}

\newcommand{\Algref}[1]{\hyperref[#1]{Algorithm~\ref*{#1}}}

\newcolumntype{M}[1]{>{\centering\arraybackslash}m{#1}}

\newcommand{\Stampede}{\href{https://www.tacc.utexas.edu/systems/stampede2}{Stampede2}}

\newcommand{\Intercepted}{\textsc{Intercepted}}
\newcommand{\Inactive}{\textsc{Inactive}}
\newcommand{\Exterior}{\textsc{Exterior}}
\newcommand{\Interior}{\textsc{Interior}}
\newcommand{\FalseIntercepted}{\textsc{FalseIntercepted}}
\newcommand{\TrueIntercepted}{\textsc{TrueIntercepted}}
\newcommand{\Neighbors}{\textsc{NeighborsFalseIntercepted}}
\definecolor{ActiveElement}{RGB}{147,194,74}
\definecolor{InterceptedElement}{RGB}{255,208,48}
\definecolor{FalseInterceptedElement}{RGB}{92,91,255}
\newcommand{\mvec}{\textsc{matvec}}

\definecolor{cpu1}{HTML}{4CAF50}
\definecolor{cpu2}{HTML}{FFC107}
\definecolor{cpu3}{HTML}{F44336}
\definecolor{cpu4}{HTML}{2196F3}
\definecolor{cpu5}{HTML}{9932CC}

\definecolor{gpu1}{HTML}{A5D6A7}
\definecolor{gpu2}{HTML}{FFE082}
\definecolor{gpu3}{HTML}{EF9A9A}
\definecolor{gpu4}{HTML}{90CAF9}

\definecolor{sq_b1}{RGB}{37,52,148}
\definecolor{sq_b2}{RGB}{44,127,184}
\definecolor{sq_b3}{RGB}{65,182,196}
\definecolor{sq_b4}{RGB}{127,205,187}
\definecolor{sq_b5}{RGB}{199,233,180}
\definecolor{sq_b6}{RGB}{255,255,204}

\definecolor{sq_r1}{RGB}{189,0,38}
\definecolor{sq_r2}{RGB}{240,59,32}
\definecolor{sq_r3}{RGB}{253,141,60}
\definecolor{sq_r4}{RGB}{254,178,76}
\definecolor{sq_r5}{RGB}{254,217,118}
\definecolor{sq_r6}{RGB}{255,255,178}

\definecolor{sq_g1}{RGB}{0,104,55}
\definecolor{sq_g2}{RGB}{49,163,84}
\definecolor{sq_g3}{RGB}{120,198,121}
\definecolor{sq_g4}{RGB}{173,221,142}
\definecolor{sq_g5}{RGB}{217,240,163}
\definecolor{sq_g6}{RGB}{255,255,204}

\definecolor{div_c1}{RGB}{230,171,2}
\definecolor{div_c2}{RGB}{102,166,30}
\definecolor{div_c3}{RGB}{231,41,138}
\definecolor{div_c4}{RGB}{117,112,179}
\definecolor{div_c5}{RGB}{217,95,2}
\definecolor{div_c6}{RGB}{27,158,119}
\definecolor{div_c7}{RGB}{215,48,39}

\definecolor{div_d1}{RGB}{215,25,28}
\definecolor{div_d2}{RGB}{253,174,97}
\definecolor{div_d3}{RGB}{255,255,191}
\definecolor{div_d4}{RGB}{171,217,233}
\definecolor{div_d5}{RGB}{44,123,182}

\definecolor{lineclr}{RGB}{0,0,0}
\definecolor{utorange}{RGB}{0,0,255}
\definecolor{utsecblue}{RGB}{255,255,0}
\definecolor{utsecgreen}{RGB}{255,0,0}
\definecolor{red!15}{RGB}{0,255,255}
\definecolor{fillclr5}{RGB}{0,255,0}
\definecolor{fillclr6}{RGB}{255,0,255}
\definecolor{fillclr7}{RGB}{255,255,255}
\definecolor{fillclr8}{RGB}{0,0,0}

\definecolor{armygreen}{rgb}{0.29, 0.33, 0.13}
\definecolor{aurometalsaurus}{rgb}{0.43, 0.5, 0.5}
\definecolor{applegreen}{rgb}{0.55, 0.71, 0.0}
\definecolor{darkgreen}{rgb}{0.0, 0.4, 0.25}

\newcommand{\tikzcircle}[2][red,fill=red]{\tikz[baseline=-0.5ex]\draw[#1,radius=#2] (0,0) circle ;}%

\def\drawcubeI(#1,#2,#3,#4,#5){ 
\coordinate (O) at (#1,#2,#3);
\coordinate (A) at (#1,#2+#4,#3);
\coordinate (B) at (#1,#2+#4,#3+#4);
\coordinate (C) at (#1,#2,#3+#4);
\coordinate (D) at (#1+#4,#2,#3);
\coordinate (E) at (#1+#4,#2+#4,#3);
\coordinate (F) at (#1+#4,#2+#4,#3+#4);
\coordinate (G) at (#1+#4,#2,#3+#4);
\draw[#5] (O) -- (C) -- (G) -- (D) -- cycle;
\draw[#5] (O) -- (A) -- (E) -- (D) -- cycle;
\draw[#5] (O) -- (A) -- (B) -- (C) -- cycle;
\draw[#5] (D) -- (E) -- (F) -- (G) -- cycle;
\draw[#5] (C) -- (B) -- (F) -- (G) -- cycle;
\draw[#5] (A) -- (B) -- (F) -- (E) -- cycle;
}

\def\drawcubeII(#1,#2,#3,#4,#5,#6,#7){ 
\coordinate (O) at (#1,#2,#3);
\coordinate (A) at (#1,#2+#4,#3);
\coordinate (B) at (#1,#2+#4,#3+#4);
\coordinate (C) at (#1,#2,#3+#4);
\coordinate (D) at (#1+#4,#2,#3);
\coordinate (E) at (#1+#4,#2+#4,#3);
\coordinate (F) at (#1+#4,#2+#4,#3+#4);
\coordinate (G) at (#1+#4,#2,#3+#4);
\draw[#5,fill=#6,opacity=#7] (O) -- (C) -- (G) -- (D) -- cycle;
\draw[#5,fill=#6,opacity=#7] (O) -- (A) -- (E) -- (D) -- cycle;
\draw[#5,fill=#6,opacity=#7] (O) -- (A) -- (B) -- (C) -- cycle;
\draw[#5,fill=#6,opacity=#7] (D) -- (E) -- (F) -- (G) -- cycle;
\draw[#5,fill=#6,opacity=#7] (C) -- (B) -- (F) -- (G) -- cycle;
\draw[#5,fill=#6,opacity=#7] (A) -- (B) -- (F) -- (E) -- cycle;
}

\def\drawNodes(#1,#2,#3,#4,#5,#6,#7){ 
\foreach \x in {#1,#7,...,#2}{
	\foreach \y in {#3,#7,...,#4}{
		\foreach \z in {#5,#7,...,#6}{
				\draw[fill=red!60] (\x,\y,\z) circle (0.15);
				}
			}
	}				
		
}

\pgfplotsset{
  log x ticks with fixed point/.style={
      xticklabel={
        \pgfkeys{/pgf/fpu=true}
        \pgfmathparse{exp(\tick)}%
        \pgfmathprintnumber[fixed relative, precision=3]{\pgfmathresult}
        \pgfkeys{/pgf/fpu=false}
      }
  },
  log y ticks with fixed point/.style={
      yticklabel={
        \pgfkeys{/pgf/fpu=true}
        \pgfmathparse{exp(\tick)}%
        \pgfmathprintnumber[fixed relative, precision=3]{\pgfmathresult}
        \pgfkeys{/pgf/fpu=false}
      }
  }
}

\makeatletter
\newcommand\resetstackedplots{
\makeatletter
\pgfplots@stacked@isfirstplottrue
\makeatother
\addplot [forget plot,draw=none] coordinates{(48,0) (96,0) (192,0) (384,0) (768,0) (1536,0) (3072,0) (6144,0)};
}
\makeatother

\makeatletter
\newcommand\resetstackedplotsOne{
\makeatletter
\pgfplots@stacked@isfirstplottrue
\makeatother
\addplot [forget plot,draw=none] coordinates{(384,0) (768,0) (1536,0) (3072,0) (6144,0)};
}
\makeatother

\makeatletter
\newcommand\resetstackedplotsTwo{
\makeatletter
\pgfplots@stacked@isfirstplottrue
\makeatother
\addplot [forget plot,draw=none] coordinates{(16,0) (32,0) (64,0) (128,0) (256,0) (512,0) (1024,0) (2048,0) (4096,0) (8192,0) (16384,0) (32768,0)};
}
\makeatother

\makeatletter
\newcommand\resetstackedplotsThree{
\makeatletter
\pgfplots@stacked@isfirstplottrue
\makeatother
\addplot [forget plot,draw=none] coordinates{(2,0) (4,0) (8,0) (16,0) (32,0) (64,0)};
}
\makeatother

\makeatletter
\newcommand\resetstackedplotsFour{
\makeatletter
\pgfplots@stacked@isfirstplottrue
\makeatother
\addplot [forget plot,draw=none] coordinates{(4,0) (8,0) (16,0) (32,0) (64,0)};
}
\makeatother

\makeatletter
\newcommand\resetstackedplotsFive{
\makeatletter
\pgfplots@stacked@isfirstplottrue
\makeatother
\addplot [forget plot,draw=none] coordinates{(1,0) (2,0) (4,0) (8,0) (16,0) (32,0) (64,0) (128,0)};
}
\makeatother

\makeatletter
\newcommand\resetstackedplotsSix{
\makeatletter
\pgfplots@stacked@isfirstplottrue
\makeatother
\addplot [forget plot,draw=none] coordinates{(2,0) (4,0) (8,0) (16,0) (32,0) (64,0) (128,0)};
}
\makeatother

\allowdisplaybreaks

\journal{Journal of Computational Physics}

\begin{document}

\begin{frontmatter}

\title{Optimal Surrogate Boundary Selection and Scalability Studies \\ for the Shifted Boundary Method on Octree Meshes}
\author[ISU]{Cheng-Hau Yang\texorpdfstring{\corref{eqb}}{}}
\ead{chenghau@iastate.edu}
\author[ISU]{Kumar Saurabh\texorpdfstring{\corref{eqb}}{}}
\ead{maksbh@iastate.edu}
\author[Duke]{Guglielmo Scovazzi}
\ead{guglielmo.scovazzi@duke.edu}
\author[poliTO]{Claudio Canuto}
\ead{claudio.canuto@polito.it}
\author[ISU]{Adarsh Krishnamurthy\texorpdfstring{\corref{cor}}{}}
\ead{adarsh@iastate.edu}
\author[ISU]{Baskar Ganapathysubramanian\texorpdfstring{\corref{cor}}{}}
\ead{baskarg@iastate.edu}
\cortext[cor]{Corresponding authors}
\cortext[eqb]{These author contributed equally}
\address[ISU]{Department of Mechanical Engineering, Iowa State University, Ames, IA}
\address[Duke]{Department of Civil and Environmental Engineering, Duke University, Durham, North Carolina 27708, USA}
\address[poliTO]{Dipartimento di Scienze Matematiche, Politecnico di Torino, Corso Duca degli Abruzzi 24, 10129 Torino, Italy}

\begin{abstract}
The accurate and efficient simulation of Partial Differential Equations (PDEs) in and around arbitrarily defined geometries is critical for many application domains. Immersed boundary methods (IBMs) alleviate the usually laborious and time-consuming process of creating body-fitted meshes around complex geometry models (described by CAD or other representations, e.g., STL, point clouds), especially when high levels of mesh adaptivity are required. In this work, we advance the field of IBM in the context of the recently developed Shifted Boundary Method (SBM). In the SBM, the location where boundary conditions are enforced is shifted from the actual boundary of the immersed object to a nearby surrogate boundary, and boundary conditions are corrected utilizing Taylor expansions. This approach allows choosing surrogate boundaries that conform to a Cartesian mesh without losing accuracy or stability. Our contributions in this work are as follows: (a) we show that the SBM numerical error can be greatly reduced by an optimal choice of the surrogate boundary, (b) we mathematically prove the optimal convergence of the SBM for this optimal choice of the surrogate boundary, (c) we deploy the SBM on massively parallel octree meshes, including algorithmic advances to handle incomplete octrees, and (d) we showcase the applicability of these approaches with a wide variety of simulations involving complex shapes, sharp corners, and different topologies. Specific emphasis is given to Poisson's equation and the linear elasticity equations. 
\end{abstract}


\end{frontmatter}


\section{Introduction}
\label{Sec:Intro}
Accurate numerical solution of PDEs in and around complex objects has a significant impact on various problems in science and technology. Examples include structural analysis of complex architectures, thermal analysis over complex geometries in semiconductor electronics, and flow analysis over complex geometries in aerodynamics. Standard numerical approaches for solving these PDEs on complex geometries--finite difference method (FDM), finite element method (FEM), or finite volume method (FVM)--usually rely on the generation of body-fitted meshes. This is a major bottleneck, as creating an analysis-suitable body-fitted mesh with appropriate refinement around the complex geometry is usually time-consuming and labor-intensive. This issue is exacerbated in problems involving moving bodies or multiphysics couplings, for which deforming meshes or re-meshing is often required (sometimes at every time step).

Immersed boundary methods (IBM) alleviate the requirement of body-fitted meshes by relaxing the requirement that the mesh conforms to the object \citep{mittal2005immersed,peskin1972flow}. IBM allowed scalable mesh generation, such as a Cartesian grid or tree-based approaches (quadtree/octree), to be deployed for simulating PDEs in and around complex objects. In this work, we concentrate on IBM in the context of FEM-based discretizations. Two main flavors of IBMs exist in this FEM context: immersogeometric analysis (IMGA, an acronym that will also refer, in what follows, to cutFEMs, the Finite Cell Method, and related approaches) and the Shifted Boundary Method (SBM). 

In IMGA, the boundary representation of the body (B-rep, NURBS, or STL) is immersed into a non-body-fitted spatial discretization. The Dirichlet boundary conditions are enforced weakly on the immersed boundary surfaces using Nitsche's method, which proved a flexible, robust and consistent approach. Interested readers are referred to \citep{xu2016tetrahedral,HOANG2019421,DEPRENTER2019604,ZhuQiming201911,saurabh2021industrial,Hsu-2016-IMGA,wang2017rapid,balu2023direct,XU2021103604,Parvizian:07.1,massing2015nitsche,burman2015cutfem,burman2010ghost,burman2014fictitious,schott2015face,schott2014new,burman2012fictitious,burman2014unfitted} for a detailed discussion of the mathematical formulation and practical deployment of the IMGA. The IMGA has been deployed to solve several industrial-scale complex problems~\cite{saurabh2021industrial}, but suffers from the following drawbacks:
\begin{itemize}[topsep=3pt,itemsep=0pt,left=0pt]
    \item \textbf{Sliver cut-cells}: The presence of sliver cut-cells (i.e., elements intersected by the object boundary that contain a very small volume of the object) may significantly deteriorate the conditioning of the algebraic system of equations. Literature suggests removing these so-formed sliver cut cells from the global assembly can prevent such deterioration in conditioning. However, this comes at the cost of accuracy. Alternatively, there have been studies demonstrating the design of preconditioners to alleviate this issue~\citep{DEPRENTER2019604}. But, this has been limited to simpler operators such as Poisson's and Stokes and requires the development of preconditioners for other PDEs. Sliver-cut cells can even produce a loss of numerical stability~\cite{burman2010ghost,burman2015cutfem}.
    
    \item \textbf{Load balancing}: The accuracy of the IMGA is strongly contingent upon the accurate integration of cut cells. Accurate integration is performed by increasing the number of quadrature points in the cut elements. However, this leads to the issue of load balancing when performing parallel (distributed memory) simulations, as different elements end up having different amounts of computations. Furthermore, it also invalidates the tensor structure of the basis function that can be exploited to optimize the matrix and vector assembly \citep{saurabh2021industrial,saurabh2022scalable}.
\end{itemize}

The Shifted Boundary Method (SBM)~\citep{main2018shifted,Main2018TheSB,KARATZAS2020113273,atallah2020second,atallah2021shifted,atallah2021analysis,colomes2021weighted,atallah2022high,ZENG2022115143} alleviates the aforementioned IMGA issues. The central idea of SBM is to impose the boundary conditions \textit{not on the true boundary} ($\G$, see Fig.~\ref{fig:SBM_domains_NA}) \textit{but rather on a surrogate boundary in proximity of the true boundary} ($\tG$, see again Fig.~\ref{fig:SBM_domains_NA}). The appropriate value of the applied boundary condition is determined by performing a Taylor series expansion. The surrogate boundary and associated shifted boundary conditions essentially transform the problem of solving the PDE in the complex original domain (denoted as $\Om$) into a body-fitted problem in the surrogate domain (denoted as $\tO$). This strategy overcomes the challenges associated with IMGA approaches. The SBM differs from the IMGA in the following aspects:
\begin{itemize}[topsep=3pt,itemsep=0pt,left=0pt]
    \item In IMGA, the volume integration is performed over $\Om$; whereas in SBM it is performed over $\tO$. Therefore, IMGA requires a classification test (to classify if a Gauss quadrature point belongs to $\Om$ or $\neg \Om$) for each Gauss quadrature point in the cut elements. In contrast, SBM does not require any such test. The integration for SBM is done over all Gauss points that belong to elements within the surrogate domain $\tO$.
    \item The integration over all Gauss points in SBM eliminates the poor conditioning of discrete operators due to the sliver cut cells arising in IMGA. 
    \item Additionally, SBM requires no adaptive quadrature for maintaining accuracy. This obviates the need for special algorithmic treatments (like weighted partitioning~\citep{saurabh2021industrial}) to ensure load balancing. In addition, the tensor nature of the basis function is retained, which can be leveraged for performance enhancement using fast vector-matrix assembly.
\end{itemize}

SBM, therefore, appears to be a promising numerical method for solving PDEs over complex domains. 
In this work, we seek to address some of the relevant questions for the practical adoption of SBM---especially for simulations over complex CAD geometry domains---that are important and yet somewhat missing from the existing literature. Specifically, we address the following:
\begin{enumerate}[topsep=0pt,itemsep=0pt]
    \item We extend the numerical analysis of SBM to cover cases when the true domain is a subset of the surrogate domain ($\Om \subset \tO$).  
    \item Different (cartesian aligned) surrogate boundaries can be constructed for a complex domain. However, some of these boundaries can be invalid, leading to disconnected surrogate domains. In this work, we codify the requirements for the set of edges/faces (in two/three dimensions) to form a valid surrogate boundary.
    \item Among these possible candidate surrogate domains, we identify the optimal surrogate domain, with boundary $\tG$, that exhibits the best accuracy. We define a simple, scalable strategy to identify this optimal surrogate boundary. 
    \item We develop the data structures and algorithms required for the scalable deployment of SBM on adaptive, incomplete octree grids. We illustrate good scaling behavior of the framework and showcase the utility of the framework by simulating a wide variety of complex three-dimensional shapes.
\end{enumerate}
The present work focuses on a particular class of PDEs, namely elliptic PDEs with applications involving diffusion (Poisson's equation) and structural mechanics (linear elasticity) problems. The remaining paper is organized as follows: In \secref{sec:Math}, we describe the mathematics of the SBM along with a description of the surrogate boundary. In \secref{sec:OptimalSurrogateBoundary} we outline the definition, approach, and algorithms for identifying the optimal surrogate boundary. In \secref{sec:AlgorithmAndImplement} we provide the details of the algorithms for the scalable deployment of SBM. In \secref{sec: Results}, we illustrate this framework with extensive numerical examples in two and three dimensions. We summarize conclusions in \secref{Sec:Conclusions}.


\section{Mathematical Formulations}
\label{sec:Math}

\subsection{The immersed variational formulation over the physical domain}
Consider the non-homogeneous elliptic equation
\begin{equation}
\label{eq:Poisson}
\begin{split}
   -\Delta u & =  f \quad \mathrm{on} \quad \Om \; , \\
    u & = u_D \quad \mathrm{on} \quad \Gamma_D \; ,
\end{split}
\end{equation}
where we are interested in solving for the scalar field, $u$, over the (immersed) domain of interest $\Om$ with boundary $\partial \Om=\Gamma_D$. Defining the appropriate functional spaces for test and trial functions, the weak formulation for Poisson's problem can be written as:
\begin{equation}
    (\nabla u_h,\nabla w_h)_{\Om_h} =  
    (f,w_h)_{\Om_h}
    + \underbrace{\avg{ \nabla u_h \cdot \bs{n} \, , \,  w_h}_{\Gamma_{D,h}}}_{\mathrm{Consistency \; term}}
    + \underbrace{\avg{u_h - u_D \, , \,   \nabla w_h \cdot \bs{n}}_{\Gamma_{D,h}}}_{\mathrm{Adjoint\; consistency \; term}} 
    - \underbrace{\avg{ \alpha \,  h^{-1} \,  (u_h - u_D) \, ,  \, w_h}_{\Gamma_{D,h}}}_{\mathrm{Penalty \; term}} ,
    \label{eq:weakIMGA}
\end{equation}
 where $\alpha$ is the penalty parameter for the Dirichlet boundary condition of the Poisson's equation, and $h$ is the element size.
 
The last three terms in \eqnref{eq:weakIMGA}---consistency, adjoint consistency, and penalty terms---
are the result of weakly applying the Dirichlet boundary condition as a surface integral. These extra terms result in surface integration over the true geometry, assuming that the finite element interpolation space can describe it exactly (otherwise a geometric discretization error may be introduced). 

In addition to the scalar elliptic equation (Poisson equation), we also consider the equations of linear elasticity. Here, we are interested in solving for the displacement vector field, $\bs{u}$. There are three essential equations for static linear elasticity. First, the equilibrium equation (and associated boundary conditions): 
\begin{equation}
\begin{split}
\nabla \cdot \bs{\sigma} + \bs{b} &= 0, \quad \mathrm{on} \quad \Om \; , \\
\bs{u} &= \bs{u}_D \quad \mathrm{on} \quad \Gamma_D \; ,
    \label{eq:equilibrium}
\end{split}
\end{equation}
where $\bs{\sigma}$ is the stress tensor, $\bs{b}$ is the body force, and again $\partial \Om = \Gamma_D$.
Second, the kinematics equation:
\begin{equation}
    \bs{\varepsilon}(\bs{u}) = \nabla^s \bs{u} = \frac{1}{2} (\nabla \bs{u} + (\nabla \bs{u})^T) \; ,
    \label{eq:kinematics}
\end{equation}
where $\bs{\varepsilon}(\bs{u})$ is the strain tensor and $\bs{u}$ is the displacement vector.
Third, the constitutive equation:
\begin{equation}
    \bs{\sigma} = \bs{C}\bs{\varepsilon}(\bs{u}) \; ,
    \label{eq:constitutive}
\end{equation}
where $\bs{C}$ is the elastic stiffness tensor. For isotropic materials, $\bs{C}$ can be written as a combination of Young's modulus $E$ and Poisson's ratio $\nu$. 
Integrating by parts and using Nitsche's method to weakly enforce the Dirichlet boundary conditions, the variational form of linear elasticity can be stated as:
\begin{equation}
    ( \bs{C} \bs{\varepsilon}(\bs{u}) \, , \, \nabla^s \bs{w}_h )_{\Om_h} =  
    (\bs{b} \, , \, \bs{w}_h )_{\Om_h} 
    + \underbrace{\avg{ (\bs{C} \bs{\varepsilon}(\bs{u})) \cdot \bs{n} \, , \,  \bs{w}_h}_{\Gamma_{D,h}}}_{\mathrm{Consistency \; term}}
    - \underbrace{\avg{ \bs{u}_h - \bs{u}_D \, , \, (\bs{C}\nabla^s \bs{w}_h) \cdot \bs{n}}_{\Gamma_{D,h}}}_{\mathrm{Adjoint\; consistency \; term}} 
    - \underbrace{\avg{ \gamma \,  h^{-1} \,  (\bs{u}_h - \bs{u}_D) \, ,  \, \bs{w}_h}_{\Gamma_{D,h}}}_{\mathrm{Penalty \; term}} \; ,
    \label{eq:weakIMGALE}
\end{equation}
where the $\gamma$ is the penalty parameter for the Dirichlet boundary condition of the linear elasticity, and $h$ is the element size. The appropriate function spaces are used for the solution $\bs{u}$ and test function $\bs{w}$.

\begin{figure}[t!]
    \centering
    \begin{subfigure}{0.25\textwidth}
        \includegraphics[width=\linewidth,trim=0 0 0 0,clip]      {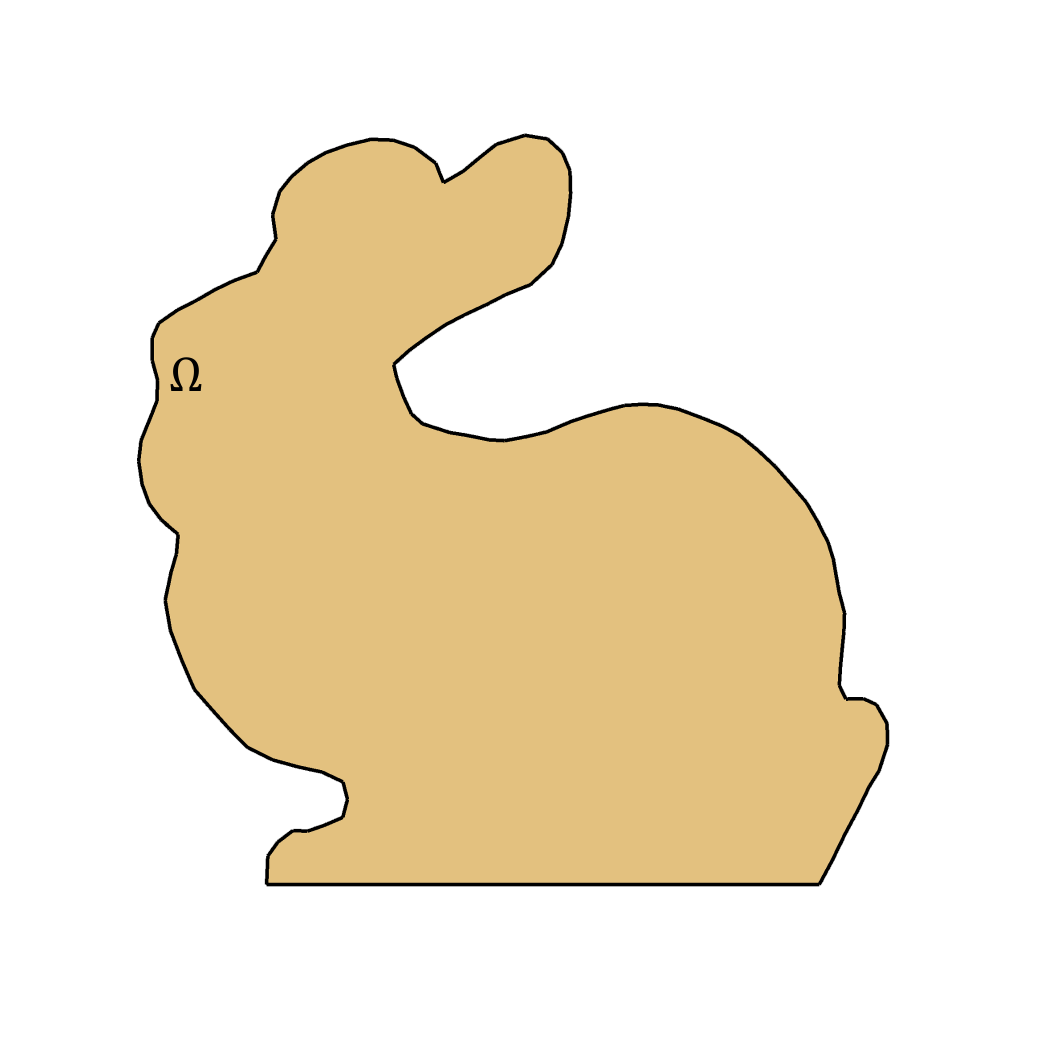}
        \caption{$\Om$ and $\G = \partial \Om$}
    \end{subfigure}
    \begin{subfigure}{0.25\textwidth}
        \includegraphics[width=\linewidth,trim=0 0 0 0,clip]      {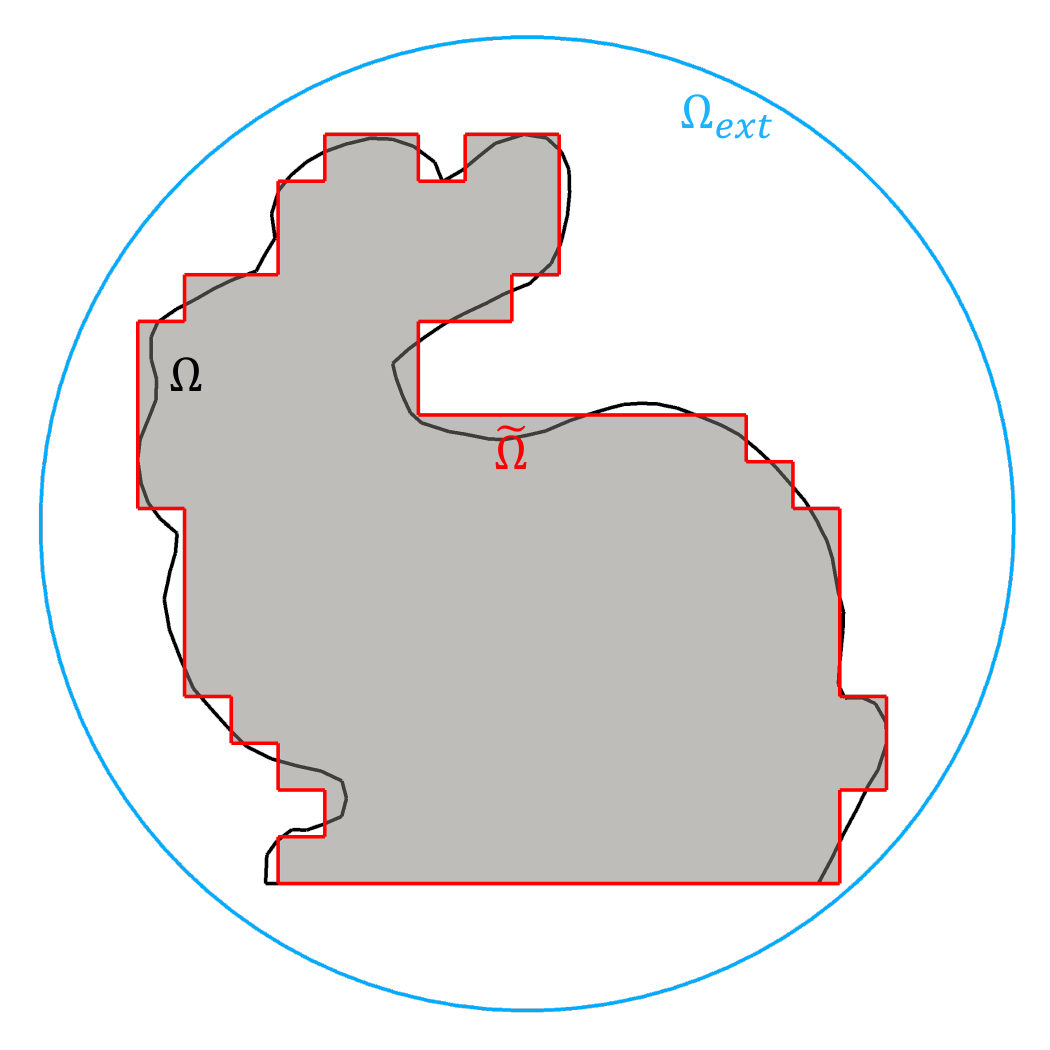}
        \caption{$\Om$, $\tO$, and $\Oe$}
    \end{subfigure}
    \caption{
    Domain definitions for the SBM numerical analysis. The domains are classified into three types, with a corresponding color scheme:
    a) The physical (or true) domain $\Om$ (\textcolor[RGB]{227,193,127}{$\blacksquare$}), enclosed by the physical (or true) boundary $\Gamma$ (\textcolor{black}{---});    
    b) The surrogate domain $\tO$ (\textcolor[RGB]{198,189,186}{$\blacksquare$}), enclosed by the surrogate boundary $\tG$ (\textcolor{red}{---}); 
    c) The extended domain $\Oe$, enclosed by the blue circle (\textcolor[RGB]{66,192,255}{---}). 
    The extended domain $\Oe \supset \Om$ will only be used in mathematical proofs. 
    } 
    \label{fig:SBM_domains_NA}
\end{figure}

\subsection{The variational formulation for the Shifted Boundary Method over the surrogate domain}
The SBM introduced in~\citep{main2018shifted} discretizes the governing equations on a surrogate domain $\tO$ of boundary $\tG$ (rather than $\Om$ and $\Gamma$), where $\tO$ and $\tG$ do not contain any cut elements or cut element sides, respectively. 
For example, looking at the sketchs in \figref{fig:SBM_domains_NA}, $\Om$ is enclosed by $\G$ (the black curve), while $\tO$ is enclosed by $\tG$ (the red segmented curve).
The SBM resorts to a Taylor expansion of the solution variable at the surrogate boundary to {\it shift} the value of the boundary condition from $\G$ to $\tG$. It is important to note that the choices of $\tO$ and $\tG$ are not independent but must satisfy certain constraints, later discussed in \secref{sec: surrogate_math}. Enforcing the Dirichlet boundary condition $u = u_D$ on $\Gamma_D$ through the SBM, we deduce the following Galerkin discretization of the Poisson equation as shown below. Here, $V_h$ represents the appropriate function space, and subscript $h$ represents the finite dimensional analogue of operators/domains after discretization with a tesselation of size $h$.

\begin{quote}
Find $u_h \in V_h(\tO)$ such that, $\forall w_h \in V_h(\tO)$
\begin{align}
\label{eq:SBM}
( \nabla u_h \, , \, \nabla w_h )_{\tO} 
& = \;
(f \, , \, w_h )_{\tO} 
 + \underbrace{\avg{ \nabla u_h  \cdot \ti{\bs{n}} \, , \,  w_h}_{\tGD}}_{\mathrm{Consistency \; term}}
 + \underbrace{\avg{ \tS u_h - u_D \, , \, \nabla w_h \cdot \ti{\bs{n}}}_{\tGD}}_{\mathrm{Adjoint\; consistency \; term}}
 \\
& \phantom{=} \; \; \;
- \underbrace{\avg{ \alpha \,  h^{-1} \,  (\tS u_h-u_D) \, , \,  \tS w_h }_{\tGD}}_{\mathrm{Penalty \; term}}
\, ,
\end{align}
where $\tS v$ is the boundary shift operator:
\begin{equation}\label{eq:def-bndS}
\tS v := v+ \nabla v \cdot \bs{d} \; , \qquad  \text{on \ } \tGD \, .
\end{equation}
 \end{quote}
 
Similarly, the SBM Galerkin discretization for static linear elasticity with Dirichlet boundary condition $\bs{u} = \bs{u}_D$ on $\Gamma_D$ can be stated as:
\begin{quote}
Find $\bs{u}_h \in V_h(\tO)$ such that, $\forall \bs{w}_h \in V_h(\tO)$
\begin{align}
\begin{split}
\label{eq:SBMLE}
 ( \bs{C} \bs{\varepsilon}(\bs{u}) \, , \, \nabla^s \bs{w}_h )_{\tO} 
& =
(\bs{b} \, , \, \bs{w}_h )_{\tO} 
+\underbrace{\avg{ (\bs{C} \bs{\varepsilon}(\bs{u})) \cdot \ti{\bs{n}} \, , \,  \bs{w}_h}_{\tGD}}_{\mathrm{Consistency \; term}}
- \underbrace{\avg{ \tSLE \bs{u}_h - \bs{u}_D  \, , \, (\bs{C}\nabla^s \bs{w}_h) \cdot \ti{\bs{n}}}_{\tGD}}_{\mathrm{Adjoint \; Consistency \; term}}
\\
&-\underbrace{\avg{ \gamma \,  h^{-1} \,  (\tSLE \bs{u}_h - \bs{u}_D)\, , \,  \tSLE \bs{w}_h }_{\tGD}}_{\mathrm{Penalty \; term}} .
\end{split}
\end{align}
 \end{quote}

The adjoint consistency and penalty terms in \eqnref{eq:SBM} and \eqnref{eq:SBMLE} are modified in the SBM formulation using Taylor expansions. We refer the interested reader to the detailed derivation of the formulation by \citet{atallah2021shifted}.

\subsection{Numerical Analysis of the Shifted Boundary Method over extended surrogate domains}

In order to identify the surrogate domain $\tO$ that leads to the most accurate results, we need to first understand the behavior of the SBM approximation when the surrogate domain extends beyond the physical domain $\Omega$, as shown, for example, in the sketch on the right of  \figref{fig:SBM_domains_NA}. As a starting point in the numerical analysis, we will need a number of definitions and assumptions. 

The true domain $\Om$ is assumed to have Lipschitz boundary $\Gamma=\partial \Omega$.
The surrogate domain $\tO$ -- in contrast with previous versions of the SBM -- is not necessarily contained in $\Om$, but may include elements that are cut by $\G$ (called intercepted elements in the sequel). Its boundary is indicated by $\tG$.

We then introduce two collections of elements: (a) the collection, $\ti{\cT}_h$, of all the elements $T$ of the grid that are contained in $\tO$; and (b) the collection, $\hat{\cT}_h$, of all the elements $T$ of the grid that are contained in $\hO$, where $\hO$ is the union of the elements cut by $\Omega$ or strictly contained in $\Omega$. Hence, $\tO \subseteq \hO$, but it is not necessarily true that $\tO \equiv \hO$. Here, $\hO$ can be thought of as the circumscribing cartesian mesh of $\Om$. We next define a domain $\Oe$ {\it with smooth boundary} and such that
$\mathrm{cl}(\hO) \subset \Oe$, where $\mathrm{cl}(\hO)$ indicates the {\it closure} of $\hO$.
Observe that $\hO$ and $\Oe$ are needed only in the mathematical analysis and are not needed in computations. For simplicity, the mathematical analysis will be developed only in the case of the Poisson problem, but conclusions similar to the ones outlined in what follows can be applied to the elasticity equations.

Consider the Poisson problem with non-homogeneous Dirichlet boundary conditions, that is, the problem of finding a $u \in H^1(\Om)$ that solves~\eqnref{eq:Poisson} for a given $f \in L^2(\Omega)$. 
We assume that either $f$ is defined directly over $\Oe$ or that we can construct a linear continuous extension operator ${E} : L^2(\Omega) \rightarrow L^2(\Oe)$ such that $Ef_{|\Omega}=f$ and $\|Ef\|_{L^2(\Oe)} \leq C \, \|f\|_{L^2(\Omega)}$, for any $ f \in L^2(\Omega)$. For example, $f$ can be extended by zero outside $\Omega$, but we use more advanced prolongation strategies in the numerical experiments.
We denote by $\bar{f}=Ef$ the extension of $f$ that we choose, and our goal is now to extend $u$ to $\bar{u}$ in $\Oe$. The following result holds, the proof of which is provided in~\ref{app:proof_ext}:
\begin{prop}
\label{prop:extension}
 There exists an extension $\bar{u}$ of $u$ in $\Oe$, such that:   
\begin{enumerate}
    \item[a)] $-\Delta \bar{u} = \bar{f}$ in $\Oe$; and 
    \item[b)] if $u \in H^2(\Omega)$, then $\bar{u} \in H^2(\tO)$, with 
    $\| \bar{u} \|_{H^2(\tO)} \leq C \; \|u\|_{H^2(\Omega)}$. 
\end{enumerate}
\end{prop}

The importance of having extensions $\bar{u}$ and $\bar{f}$ of $u$ and $f$ that satisfy conditions a) and b) above is needed when studying the convergence of the SBM for a surrogate domain $\tO$ that is not completely contained in the physical domain $\Omega$. 
Observe that the numerical stability of the SBM is not affected by the particular choice of surrogate domain, as long as $\bs{d}$ goes to zero as the grid size $h$ is refined.
We state then the following result without proof, since the derivations will not differ from the ones already found in the existing literature on SBM~\cite{atallah2020second}.
\begin{thm}[Coercivity]
\label{thm:coerc-aPoisson}
Consider the bilinear form $a_h(u_h\, , \, w_h)$ defined in~\eqref{eq:SB_Poisson_uns}
and assume there exist constants $c_d>0$ and $\zeta >0$ such that
	\begin{equation}\label{eq:ass_lbb}
	\| \, \bs{d}(\ti{\bs{x}}) \, \| \leq c_d \, h_T \, \hat{h}_T^{\zeta} \;  \qquad \forall \ti{\bs{x}} \in \tG \cap T, \ \ T \in \ti{\cT}_h \; , 
	\end{equation}
	where
	\begin{equation}\label{eq:hhat}
	\hat{h}_T =  l(\tO)^{-1} \,  h_T   \; .
	\end{equation}
Then, if the parameter $\alpha$ is sufficiently large and $\hat{h}_{\tG}$ sufficiently small, there exists a constant $C_a>0$ independent of the mesh size, such that	
\begin{equation} \label{eq:17}
a_h(u_h\,, \,u_h) \geq  C_{a} \, \| \, u_h \, \|^2_{a}  \qquad \forall u_h \in V_h(\tO) \; ,
\end{equation}
where $ \| \, u_h \, \|_{a}^2  = \| \, \nabla u_h \, \|^2_{L^2(\tO)}   + \| \, h^{-1/2} \, \tS u_h \, \|^2_{L^2(\tG)}  \; . $
\end{thm}

\begin{remark}
Condition~\eqref{eq:ass_lbb} is just a technical condition for the proofs. In fact, we take $\| \, \bs{d}(\ti{\bs{x}}) \, \| \sim \, h_T$ in all computations presented in this work, that is, mesh refinement is obtained by just subdividing every edge of the discretization into two equal-size sub-edges.
\end{remark}

The convergence analysis, which follows the same general strategy developed in~\cite{atallah2020second,atallah2021analysis,atallah2022high}, needs to be considered with more care. In particular a convergence proof is achieved using Strang's lemma, which in turn requires a result of asymptotic consistency of the SBM. Because most of the derivations are substantially similar to the ones in~\cite{atallah2020second,atallah2021analysis,atallah2022high}, we will only focus on the differences, notably the {\it asymptotic consistency estimate}.
In the present case, recasting~\eqnref{eq:SBM} as
\begin{subequations}
\begin{align}
\label{eq:SB_Poisson_uns}
    a_h(u_h,w_h) 
    &=\;
    \ell_h(w_h) \; ,
\end{align}
with
\label{eq:consistency}    
\begin{align}
    a_h(u_h,w_h) 
    &=\;
    ( \nabla u_h \, , \, \nabla w_h )_{\tO} 
    -\avg{ \nabla u_h  \cdot \ti{\bs{n}} \, , \,  w_h}_{\tG}
    -\avg{ \tS u_h \, , \, \nabla w_h \cdot \ti{\bs{n}}}_{\tG} 
    +\avg{ \alpha \,  h^{-1} \,  \tS u_h \, , \,  \tS w_h }_{\tG}
    \; ,
    \\
    \ell_h(w_h)
    &=\;
    (f \, , \, w_h )_{\tO} 
    -\avg{ u_D \, , \, \nabla w_h \cdot \ti{\bs{n}}}_{\tG} 
    +\avg{ \alpha \,  h^{-1} \,  u_D \, , \,  \tS w_h }_{\tG}
    \; ,
  \end{align}
\end{subequations}
and replacing in~\eqnref{eq:consistency} $f$ with the extension $\bar{f}$ and $u_h$ with the extension $\bar{u}$ of the exact solution $u$, we have:
\begin{align*}
    a_h(\bar{u},w_h) - \ell_h(w_h) 
    &=\;
    ( \underbrace{- \Delta \bar{u} + \bar{f}}_{\equiv 0} \, , \, \nabla w_h )_{\tO} 
    -\avg{ \underbrace{\tS \bar{u}-u_D}_{R_h \bar{u}} \, , \, \nabla w_h \cdot \ti{\bs{n}} - \alpha \,  h^{-1} \,  \tS w_h }_{\tG}
    \; ,
\end{align*}
where $R_h \bar{u}=\tS \bar{u}-u_D$ denotes the {\it residual} of the Taylor expansion.
From this, using appropriate trace inequalities, we deduce
\begin{align*}
    | a_h(\bar{u},w_h) - \ell_h(w_h) |
    & \leq \;
    C \, 
    \| R_h \bar{u} \|_{L^2(\tO)}
    \,\| w_h \|_{V(\tO;\ti{\cT}_h)}
    \; ,
\end{align*}
where
\begin{align}
\| \, v \, \|^2_{V(\tO;\ti{\cT}_h)} 
& = 
\| \, v \, \|_{a}^2
+  |\, h\, v \, |^2_{H^2(\tO;\ti{\cT}_h)}
\end{align}
is the norm associated with the infinite dimensional space 
\begin{align}
V(\tO;\ti{\cT}_h) & = V_h(\tO) +  H^2(\tO)\;  \subset H^2(\tO;\ti{\cT}_h)\, \;  . 
\end{align}
Here $V(\tO;\ti{\cT}_h)$ is an extension of the finite dimensional space $V_h(\tO)$ of globally continuous, piecewise-linear polynomials, which contains the extension of the exact solution $\bar{u}$, that is $\bar{u} \in V(\tO;\ti{\cT}_h)$. 
Here $H^2(\tO;\ti{\cT}_h)= \prod_{T \in \ti{\cT}_h} H^2(T)$ with `broken' norm $\| \cdot \|_{H^2(\tO;\ti{\cT}_h)}= \sum_{T \in \ti{\cT}_h} \| \cdot \|_{H^2(T)}$
and `broken' seminorm $| \cdot |_{H^2(\tO;\ti{\cT}_h)}= \sum_{T \in \ti{\cT}_h} | \cdot |_{H^2(T)}$.
It is easily checked that the form $a_h(\cdot,\cdot)$ is well-defined also on the space $V(\tO;\ti{\cT}_h) \times V_h(\tO)$. 

Since $\bar{u}$ has regularity $H^2$ around $\tG$, the norm of the reminder $R_h \bar{u}$ can be estimated as in the standard case in which $\tO \subset \Omega$ (see, e.g.,~\cite{atallah2020second}), leading to:

\begin{thm}[Optimal convergence in the natural norm]
\label{thm:PoissonConvergenceNatural}
Assume that $\G$ is of class $\mathcal{C}^{2}$,  $f \in L^2(\Om)$ and $u_D \in H^{3/2}(\G)$. Under the assumption of Theorem~\ref{thm:coerc-aPoisson}, and the condition that $h_{\tG}$ is sufficiently small, the numerical solution $u_h$ satisfies the following error estimate:
\begin{align}
\label{eq:convergenceEstimate}
\| \, \bar{u} - u_h \, \|_{V(\tO;\tT)}  \leq C\, h_{\tO} \, \| \, \nabla (\nabla \bar{u}) \,  \|_{L^2(\Omega)}  
	\; ,
\end{align}
where $C>0$ is a constant independent of the mesh size and the solution. 
\end{thm}

In addition, duality estimates can be derived to show that the $L^2(\tO)$-error of the discrete solution converges with rate $3/2$, which suboptimal by an order $1/2$. Note however that optimal $L^2$-error convergence rates have been observed in all computations performed to date with the SBM, for a variety of problems and differential operators. This might indicate that the available $L^2$-error estimates are not sharp.  

\clearpage 

\section{Optimal Surrogate Boundary}
\label{sec:OptimalSurrogateBoundary}

As already discussed at length, the key aspect of the SBM is the correction of the boundary conditions on the surrogate boundary, obtained by performing a Taylor series expansion.
Previous literature \citep{main2018shifted,saurabh2021scalable} has shown that convergence in the $L^2$-norm reduces to first order, when using linear basis functions without this correction, for instance. In this section, we answer the question of constructing the \emph{optimal} surrogate boundary, which gives minimal error while retaining all the expected properties of SBM. 

It is rather intuitive to recognize that the surrogate boundary with minimum distance $\bs{d}$ (in some sense) from the true boundary should be optimal. Using a wide variety of canonical examples exhibiting complex shapes and topology, we show that solving the PDE using an optimal surrogate boundary (i.e., with minimal $\bs{d}$) can produce significantly more accurate solutions compared to a non-optimal surrogate. Given a background adaptive Cartesian mesh (octree or quadtree), identification of the optimal surrogate can be stated as an optimization problem, where the goal is to minimize the distance between the true boundary and the surrogate boundary (\eqnref{eq: optim}):
\begin{equation}
        \mathrm{arg min}_{\tG} ||{\Gamma - \tG}|| 
        :=
        \mathrm{arg min}_{\tG} \int_{\tG}
        | \bs{d} \cdot \ti{\bs{n}} | \, 
        \mathrm{d}\tG \; ,
    \label{eq: optim}
\end{equation}
which corresponds to the measure of the gap between $\G$ and $\tG$.
Performing a global optimization on the surrogate boundary is a non-trivial task. We recast this global optimization into a set of element-level optimization as:
\begin{equation}
    \begin{split}
    \mathrm{arg min}_{\tG} \int_{\tG}
        | \bs{d} \cdot \ti{\bs{n}} | \, 
        \mathrm{d}\tG
        = \mathrm{arg min}_{\tG}  \left( \sum_{T \in \ti{\mathcal{T}}_h}
        \int_{\partial T \cap \tG}
        | \bs{d} \cdot \ti{\bs{n}} | \, 
        \mathrm{d}\tG
        \right)
    \; ,
    \end{split}
\end{equation}
where we recall that $\ti{\mathcal{T}}_h$ is the collection of elements in $\tO$. Converting the global optimization into a set of element-level optimizations is algorithmically useful, both from a complexity standpoint and a communication/data structure standpoint. However, performing the optimization at the local elemental level does not guarantee the satisfaction of constraints of the surrogate boundary (described in detail in \secref{sec: surrogate_math}). 
To alleviate this issue, we modify the problem represented by \eqnref{eq: optim}:
Instead of asking the question \textit{``how close is the surrogate boundary to the true boundary?"} we ask the question \textit{``how close is the surrogate volume to the true volume?"} Basically we approximate \eqnref{eq: optim} as:
\begin{equation}\label{eq:surrogate_opt_eqn}
    \begin{split}
    \mathrm{arg min}_{\tG} \int_{\tG}
        | \bs{d} \cdot \ti{\bs{n}} | \, 
        \mathrm{d}\tG
        \approx
        \mathrm{arg min}_{\tO}
        | (\Om \setminus \tO ) \cup (\tO \setminus \Om) |
    \; .
    \end{split}
\end{equation}

\subsection{Algorithmic description of the surrogate domain and its boundary}
\label{sec: surrogate_math}

\begin{figure}[t!]
    \centering
    \includegraphics[width=0.5\linewidth,trim=0 240 0 170,clip]{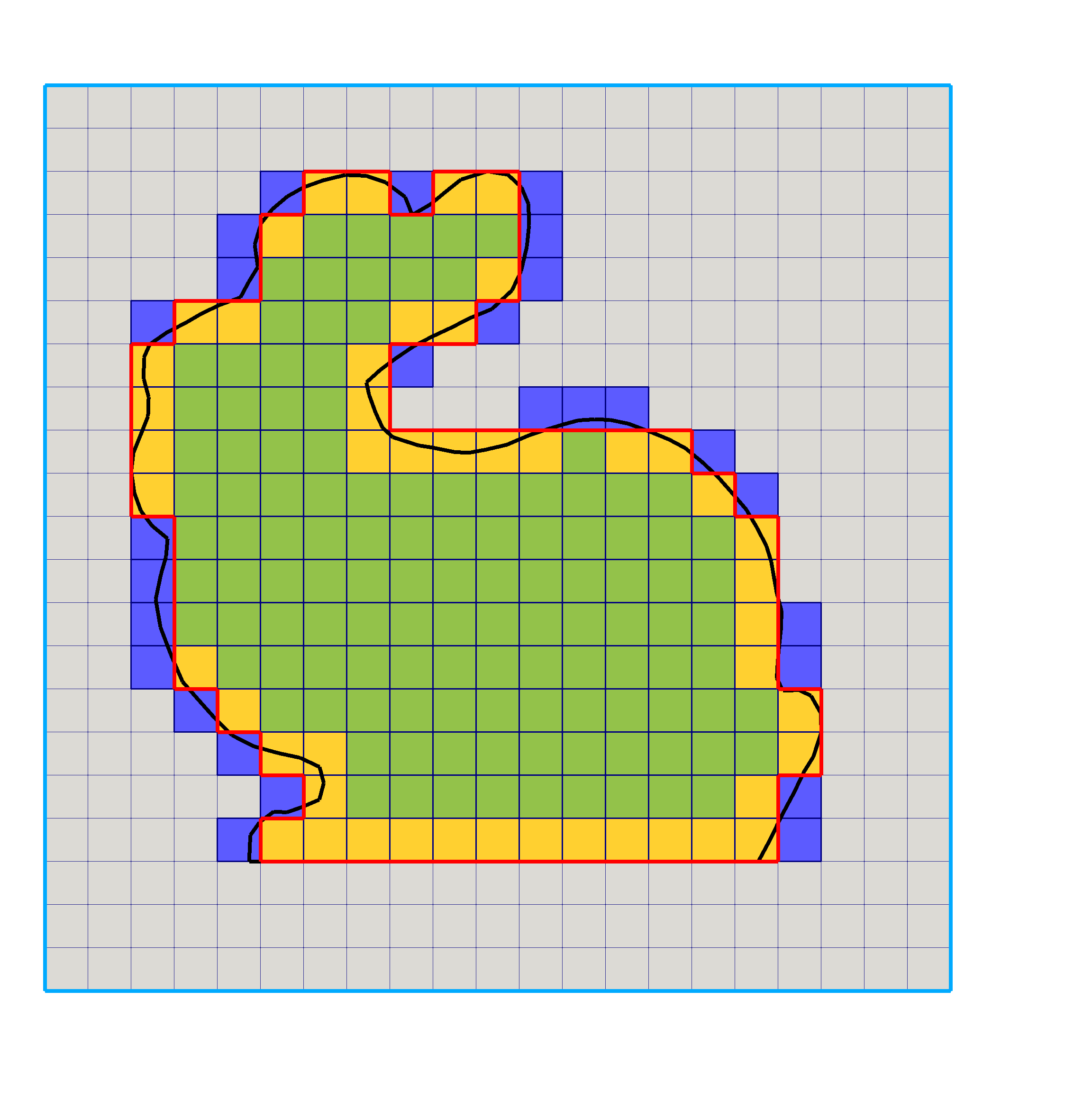}
    \captionsetup{singlelinecheck=off}
    \caption[]{
    Color scheme depicting four different types of elements in the SBM algorithms and the four types of associated domains: \Interior{} elements (\textcolor{cpu1}{$\blacksquare$}), \TrueIntercepted{} elements (\textcolor{cpu2}{$\blacksquare$}), \FalseIntercepted{} elements (\textcolor{FalseInterceptedElement}{$\blacksquare$}), and \Exterior{} elements (\textcolor[RGB]{220,218,213}{$\blacksquare$}).
    }
    \label{fig:SBM_elements_domains}
\end{figure}
\citet{main2018shifted} proposed to define the surrogate boundary as the closest projection of the true boundary. However, the surrogate domain $(\tO)$ was constructed using only elements that are completely contained in the true domain. We formulate the requirements of the surrogate domain $(\tO)$ and surrogate boundary $(\tG)$ more formally, and limit our discussion to quadrilateral/hexahedral elements. This is motivated by the fact that scalable adaptive algorithms exist for creating quad/octree meshes \citep{burstedde2011p4est,sundar2008bottom,ishii2019solving,dealII90}. 

We start with some terminology that we will use throughout the manuscript. We encourage the reader to familiarize with \figref{fig:SBM_elements_domains} before moving on to the definitions below. The figure illustrates all the domains defined in earlier sections (true, surrogate, circumscribing, extension) and relates them to the corresponding mesh elements. 
Namely:
    \begin{itemize}
    \setlength\itemsep{0em}
    \item \Interior{} elements (\textcolor{cpu1}{$\blacksquare$}): the elements whose 
    nodes are inside the physical boundary ($\Gamma$).
    \item \Exterior{} elements (\textcolor[RGB]{220,218,213}{$\blacksquare$}): the elements whose four nodal points are outside the physical boundary ($\Gamma$).
    \item \Intercepted{} elements (\textcolor{cpu2}{$\blacksquare$} $\cup$ \textcolor{FalseInterceptedElement}{$\blacksquare$}): elements whose nodal points are partially within and partially outside of the physical boundary $\Gamma$. We further subdivide \Intercepted{} elements into two categories -- \TrueIntercepted{}, and \FalseIntercepted{} -- based on whether they fall within the surrogate boundary $\tilde{\Gamma}$ (\textcolor{red}{---}).
    A strategy for this classification is provided in this section.
    The identification of the optimal surrogate boils down to determining which \Intercepted{} elements belong to each of these sub-categories:
    \begin{itemize}
    \setlength\itemsep{0em}
        \item \TrueIntercepted{} elements (\textcolor{cpu2}{$\blacksquare$}): the \Intercepted{} elements that are inside the surrogate boundary $\tilde{\Gamma}$. These elements are part of the SBM calculation.
        \item \FalseIntercepted{} elements (\textcolor{FalseInterceptedElement}{$\blacksquare$}): the \Intercepted{} elements that are outside the surrogate boundary $\tilde{\Gamma}$. These elements are not part of the SBM calculation.
    \end{itemize}
    \end{itemize}    %
    The sketch also shows the three different domains considered in what follows: 
    \begin{itemize}
    \setlength\itemsep{0em}
        \item The physical (or true) domain $\Omega$: the domain enclosed by the physical (or true) boundary $\Gamma$ (\textcolor{black}{---}). 
        \item The surrogate domain $\tilde{\Omega} = $ \textcolor{cpu1}{$\blacksquare$} $\cup$ \textcolor{cpu2}{$\blacksquare$}: the domain enclosed by the surrogate boundary $\tilde{\Gamma}$ (\textcolor{red}{---}), that is the union of the \Interior{} elements (\textcolor{cpu1}{$\blacksquare$}) and the \TrueIntercepted{} elements (\textcolor{cpu2}{$\blacksquare$}). This is the domain over which the SBM calculations are performed. 
        \item The extended domain $\Omega_{ext} $: the domain enclosed by the blue square (\textcolor[RGB]{66,192,255}{---}). This domain contains \Interior{} elements (\textcolor{cpu1}{$\blacksquare$}), \TrueIntercepted{} elements (\textcolor{cpu2}{$\blacksquare$}), \FalseIntercepted{} elements (\textcolor{FalseInterceptedElement}{$\blacksquare$}), and \Exterior{} elements (\textcolor[RGB]{220,218,213}{$\blacksquare$}).
    \end{itemize}
    
We refer the reader to \secref{sec:AlgorithmAndImplement}, which contains algorithmic details of how to perform the classification of the various element types. We next state formal definitions that allow rigorous algorithmic developments of scalable strategies for constructing these optimal surrogate domains: 

\begin{definition}[\textbf{Node}]
    A node is defined as a point $\Vec{x} \in \mathcal{R}^{dim}$, where $dim$ is the domain dimensionality (2, or 3).
\end{definition}

\begin{definition}[\textbf{Node classification}]
    A node is classified as \Interior{} if it lies within the true domain $\Om$, otherwise is classified as \Exterior{}
    \label{defn:nodeClassification}
\end{definition}

\begin{definition}[\textbf{Element node relation}]
    Each octant or element in the mesh comprises a certain number of nodes. The actual number of nodes that comprise an element depends on the order of the basis function and the dimension and varies as $(p + 1)^{dim}$, where $p$ is the basis function order, and $dim$ is the dimensionality. 
    \label{defn:ElementNodeRelation}
\end{definition}

\begin{definition}[\textbf{Element classification}]
   The elements/octants of the octree are categorized into three categories: \Interior{}, \Exterior{}, \Intercepted{}.  An element is classified as \Exterior{} if all the nodes of the element are classified as \Exterior{}. Similarly, the element is classified as \Interior{} if all the nodes of the element are \Interior{}. When the nodes of the element have some nodes classified as \Interior{} and some as \Exterior{}, the elements are classified as \Intercepted{}. 
   \label{defn:elementClassification}
\end{definition}

\begin{remark}
    The classification of \Interior{} and  \Exterior{} regions depend on the domain of interest for the PDEs. It can be the inside or outside of the enclosed geometry. For instance, if one is interested in the effect of inclusions/voids then the domain of interest is the outside of the geometry defining these voids. 
\end{remark}

In practical terms, Eqn.~\ref{eq:surrogate_opt_eqn} boils down to looping over all \Intercepted{} elements and deciding whether that \Intercepted{} element should be retained in the surrogate domain ($\tO$). A simple and effective strategy is to retain an \Intercepted{} element if it encloses enough of the true domain. To formalize this, we define an additional classification of an element, \FalseIntercepted{}. 

\begin{definition}[\textbf{\FalseIntercepted{}}]
An \Intercepted{} element is classified as \FalseIntercepted{} if the ratio of the element volume exterior to $\Om$ to the total element volume is greater (>) than the threshold factor $\lambda$. 
\end{definition}

\begin{figure}[t!]
    \centering
    \begin{subfigure}{0.23\textwidth}
        \includegraphics[width=\linewidth,trim=0 0 0 0,clip]      {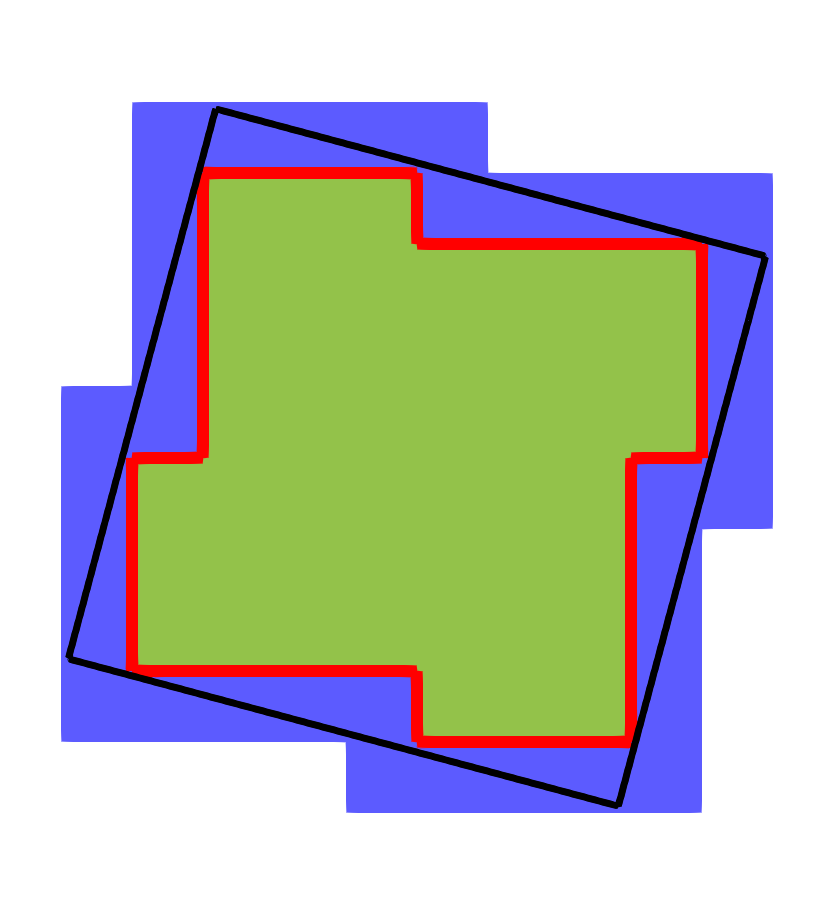}
        \caption{$\lambda$ = 0}
    \end{subfigure}
    \begin{subfigure}{0.23\textwidth}
        \includegraphics[width=\linewidth,trim=0 0 0 0,clip]      {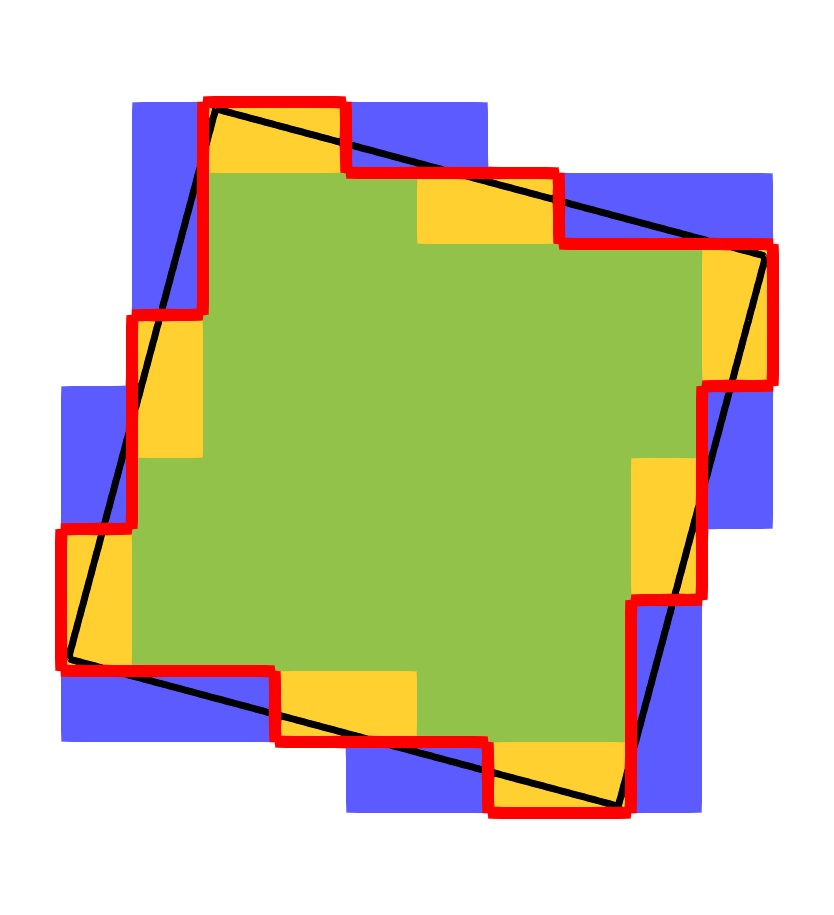}
        \caption{$\lambda$ = 0.5}
    \end{subfigure}
        \begin{subfigure}{0.23\textwidth}
        \includegraphics[width=\linewidth,trim=0 0 0 0,clip]      {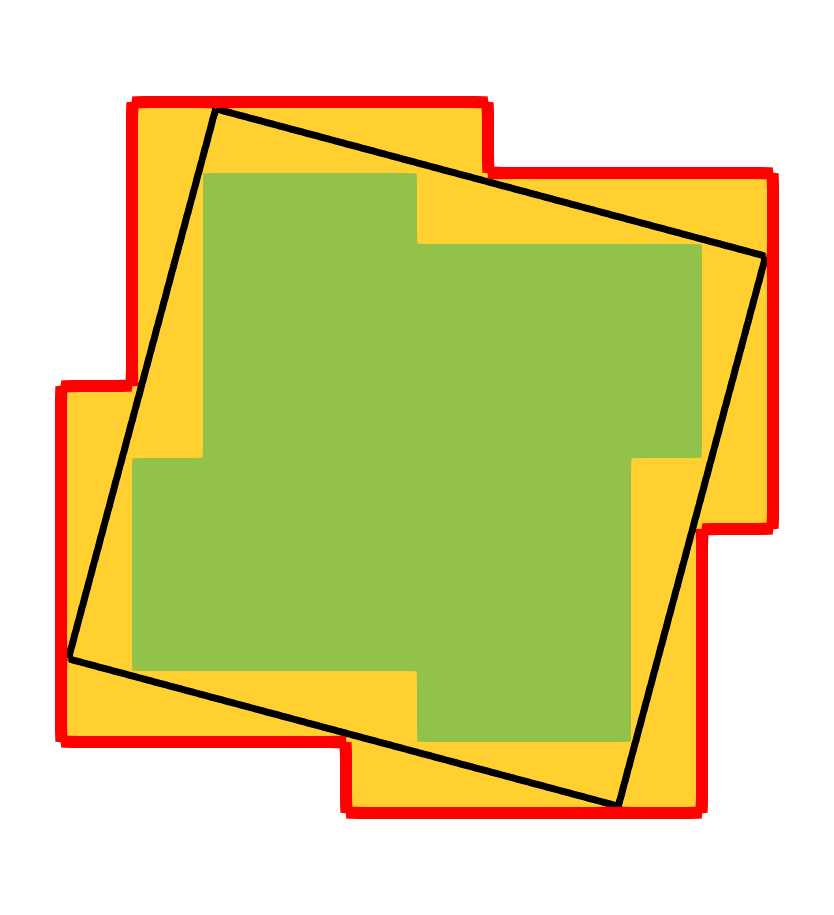}
        \caption{$\lambda$ = 1}
    \end{subfigure}
    \caption{\textbf{Surrogate boundary and the marker with varying $\lambda$}: Figure showing the classification of a) \Interior{} (\textcolor{ActiveElement}{$\blacksquare$}); b) \Intercepted{} (\textcolor{InterceptedElement}{$\blacksquare$}); c) \FalseIntercepted{} (\textcolor{FalseInterceptedElement}{$\blacksquare$}) with varying $\lambda$. \Inactive{} elements do not form the incomplete octree and are not present.  Note that (\textcolor{red}{---}) indicates the surrogate boundary whereas (\textcolor{black}{----}) denotes the true boundary.}
    \label{fig:DifferenceLambda}
\end{figure}

We note that the classification of the element as \FalseIntercepted{} is contingent on the choice of the user-defined parameter $\lambda$. When we choose $\lambda = 0$, all the \Intercepted{} elements are classified as \FalseIntercepted{}, which produces a surrogate domain that fully inscribes (i.e., is inside) the true domain. On the other hand, choosing $\lambda = 1$ leads to the inclusion of all the \Intercepted{} elements producing a surrogate domain that fully circumscribes the true domain. \figref{fig:DifferenceLambda} illustrates various surrogate boundaries as a function of varying $\lambda$. Intuitively, $\lambda=0.5$ produces an optimal surrogate that minimizes \eqnref{eq: optim}. 

The \textit{surrogate domain} $\tO$ of size $h$ is defined as a set of elements with element size $||\Delta\bs{x}|| \leq h$ such that when any extra element (of size $||\Delta\bs{x}|| \leq h$) that belongs to the complement $\tO^c$ of $\tO$ is added, {\color{black} it must be classified as either \Exterior{} or \FalseIntercepted{}}. A surrogate domain can be constructed as circumscribing ($\lambda = 1$) or inscribing ($\lambda = 0$) the true domain, or ``something in between'' ($\lambda \in (0,1)$) these two extreme cases. A \textit{surrogate boundary} is the set of faces/edges that traverses the surrogate domain $\tO$. \figref{fig:DifferenceLambda} illustrates a variety of surrogate domains and associated surrogate boundaries for a given geometry. We design algorithms such that the surrogate domain satisfies the following conditions to ensure correct computations:
\begin{itemize}
\item \textbf{Watertightness}: Ideally, {\color{black} the true boundary} must be watertight or 2-manifold as nothing can enter or leave the domain. In practice, however, the SBM approach is robust to small gaps/overlaps.
\item \textbf{Single-cycle condition}: The set of edges or faces that form the surrogate boundary must form one and exactly one cycle that traverses the surrogate domain $\tO$. In other words, there should not be any self-intersections in the surrogate boundary.
\end{itemize}

In the next section, we describe the algorithms to construct the surrogate boundary for arbitrary choices of $\lambda$ in a massively parallel environment.

\section{Algorithms and Implementation Details}
\label{sec:AlgorithmAndImplement}

\subsection{Algorithms}
To start the discussion of the algorithms for the efficient and accurate construction of surrogate domain and surrogate boundary for the SBM computations, we clarify some assumptions and motivations behind the choices described in what follows.

\subsubsection{Assumptions regarding meshes}
\label{sec: Assumption}
Before proceeding to the algorithm sections, we make the following assumption regarding the data structures.
\begin{enumerate}[itemsep=0pt,topsep=0pt]
    \item \textbf{No neighbor information:} We assume that the mesh elements do not have access neighborhood information. Tagging neighbors is particularly challenging with unstructured meshes, as elements can have varying neighbors with no plausible upper limits.
    \item \textbf{Partitioned from get-go :} The octree-based mesh data structure is partitioned right from the construction stage using distributed memory parallelism. This aspect has made octrees possible to scale to thousands of processors. This is in contrast to the traditional unstructured mesh generation, where the mesh is first generated on a single processor and later partitioned through a graph partitioning library such as \textsc{ParMetis}. This is an important aspect to consider while developing algorithms that retain the scalability of octree meshes. 
    \item \textbf{Massively parallel environment:} The algorithm proposed should scale to thousands of processors. We are not only interested in the accurate solution of PDEs but also in an efficient and scalable solution. 
    \item \textbf{Different element sizes:} Octrees can have different element sizes. We consider 2:1 balanced octrees during our algorithmic development~\cite{SundarSampathBiros08}. Additionally, we assume that the \Intercepted{} elements and \Interior{} elements that are neighbors of the \Intercepted{} elements (elements that share at least one node of the \Intercepted{} elements) are at the same level. This is done to ensure that there are no hanging nodes, which retains the simplicity of algorithms without too much extra computational cost. 
\end{enumerate}

\subsubsection{Algorithm for determining surrogate boundary for arbitrary boundary}

With the above assumption in \secref{sec: Assumption}, we can define the algorithms for determining the surrogate boundary for any arbitrary choice of $\lambda$. The basic idea of the proposed algorithm is to rely on the connection between the elements through the nodes. We note that the nodes are shared across the elements in Continuous Galerkin (CG) Finite element method. Other researchers have leveraged this to implement several graph-based algorithms for unstructured meshes, even without any neighbor information stored in the mesh data structure \citep{bogle2019parallel}. This can be efficiently performed as a series of  \mvec{} operations---a key component in FEM libraries and can be performed in a highly efficient and scalable manner~\citep{ishii2019solving,saurabh2021scalable,fernando2017machine,burstedde2011p4est}. 

\begin{algorithm}[b!]
  \footnotesize
    \caption{\textsc{SurrogateBoundaryIdentification:} Identify the surrogate boundary}
    \label{alg:surrogateBoundary}
    \begin{algorithmic}[1]
\Require Incomplete octree mesh $\mathcal{O}$,Threshold factor $\lambda$
\Ensure Surrogate boundary ($\Tilde{\Gamma}$)
\item[]
\State Marker $\mathcal{M}_1 \leftarrow $ \textsc{GenerateMarkers}($\mathcal{O}, \lambda$) \Comment{\Algref{alg:markers}}. 
\State Marker $\mathcal{M}_2$, \FalseIntercepted{} \textsc{Nodes}  $(\mathcal{N}) \leftarrow$ \textsc{GenerateNeighborsOfFalseIntercepted}($\mathcal{O}$) \Comment{\Algref{alg:neighbors}}
\State Boundary $\Tilde{\Gamma}$, Marker $M_3 \leftarrow$ \textsc{GetBoundary}($\mathcal{O},\mathcal{M}_2$) \Comment{\Algref{alg:boundary}}
\item[]
\Return ($\Tilde{{\Gamma}}$)
\end{algorithmic}
\end{algorithm}

\Algref{alg:surrogateBoundary} briefs the major step required to identify the surrogate boundary. We begin with identifying markers for each element (\Algref{alg:markers}). At this stage, each element is classified as \Interior{}, \Exterior{}, or \Intercepted{}. Next, each \Intercepted{} element is classified as \FalseIntercepted{} depending on the value of $\lambda$. For accurate evaluation of the volume term within $\Om$, we use $5^{dim}$ Gauss - Legendre points; i.e., each element is filled with 5 Gauss points in each dimension. This computational choice works well for all our results but can be easily changed at compile time.

\begin{algorithm}[t!]
  \footnotesize
    \caption{\textsc{GenerateMarkers:} Generate Marker classification (\Exterior{}, \Interior, \Intercepted, \FalseIntercepted}
    \label{alg:markers}
    \begin{algorithmic}[1]
\Require Octree $\mathcal{O}$,Threshold factor $\lambda$
\Ensure Marker $M$
\item[]
\State $\mathcal{M} \leftarrow []$
\For{element $\in \mathcal{O}$ } \Comment{Loop over the elements of octree}
\State count $\leftarrow$ 0
\For{nodes $\in$ elements} \Comment{Loop over the nodes of each element}
\If{nodes == \Interior{}} \Comment{Classify nodes as \Exterior{} or \Interior}
\State count++ \Comment{Increment for \Interior{} nodes of element}
\EndIf
\EndFor
\If{count == num\_nodes}
\State $\mathcal{M}$[element] $\leftarrow$ \Interior{}
\ElsIf{count == 0}
\State $\mathcal{M}$[element] $\leftarrow$ \Exterior{}
\Else 
\State $\mathcal{M}$[element] $\leftarrow$ \Intercepted
\State count\_gp $\leftarrow$ 0  
\Comment{Counter for number of Gauss points that are \Interior{}}
\For{gp $\in$ GaussPoints}
\If {gp $\leftarrow$ \Interior{}} 
\State count\_gp ++ 
\EndIf
\State $\lambda_c \leftarrow $ count\_gp/num\_gp \Comment{Fraction of Gauss point that are \Interior{}}
\If{$\lambda_c \geq \lambda$} \Comment{Classify elements on the basis of the threshold $\lambda$}
\State $\mathcal{M}$[element] $\leftarrow$ \FalseIntercepted
\EndIf
\EndFor
\EndIf
\EndFor
\item[]
\Return ($\mathcal{M}$)
\end{algorithmic}
\end{algorithm}

\begin{algorithm}[t!]
  \footnotesize
    \caption{\textsc{GenerateNeighborsOfFalseIntercepted:} Generate neighbors of \FalseIntercepted}
    \label{alg:neighbors}
    \begin{algorithmic}[1]
\Require Octree $\mathcal{O}$, Element marker $\mathcal{M}$
\Ensure Marker $M$ with element marker for \Neighbors, Nodal False Intercepted nodes $\mathcal{N}$
\item[]
\State $\mathcal{N}$\_ghosted $\leftarrow$ [0] \Comment{Vector of zeros with size of ghosted nodal vectors}

\For{element $\in \mathcal{O}$}
\If {element == \FalseIntercepted{}}
\For{nodes $\in$ element} 
\State $\mathcal{N}$\_ghosted[nodes] $\leftarrow$ 1 \Comment{Assign value of 1 to all the nodes of the element tagged as \FalseIntercepted}
\EndFor
\EndIf
\EndFor
\State  $\mathcal{N}$ $\leftarrow$ GhostWrite($\mathcal{N}$\_ghosted) \Comment{Ghost write}
\State  $\mathcal{N}$\_ghosted $\leftarrow$ GhostRead($\mathcal{N}$) \Comment{Ghost read}
\For{element $\in \mathcal{O}$}
\If{\texttt{anyof}(nodes) $\in$ element == 1} \Comment{If any of the nodes is marked as 1}
\State $\mathcal{M}$[element] $\leftarrow$ \Neighbors
\EndIf
\EndFor
\item[]
\Return ($\mathcal{M}$,$\mathcal{N}$)
\end{algorithmic}
\end{algorithm}

Removing \FalseIntercepted{} elements from the domain requires a change of surrogate boundary. To identify the surrogate boundary, we generate the markers for neighbors of \FalseIntercepted{}. Without the neighbor information within the mesh data structure, we rely on the efficient ~\mvec{} computation to achieve this task (\Algref{alg:neighbors}). This step can also be considered a \emph{scatter-to-gather} transformation. Each \FalseIntercepted{} element scatter the information by assigning the value of 1 to the incident nodes on that given element. In the second pass, each element gathers the data by looking into the values of the nodes that are incident on it. Once we have the nodal and elemental information about the \Neighbors{}, we can compute the faces that form the new $\tG$. Two distinct cases exist:
\begin{enumerate}[itemsep=0pt]
    \item If an element is marked as \Intercepted{}, we proceed as usual. The face(s) of the element with all the nodes marked as \Exterior{} is added to the surrogate boundary faces.
    \item For the element marked \Neighbors{}, we loop through the faces of the element. If all the nodes on a given face are either \FalseIntercepted{} or \Exterior{}, then and only then, they form the part of~$\tG$.
\end{enumerate}

\begin{algorithm}[t!]
  \footnotesize
    \caption{\textsc{GetBoundary:} Get surrogate boundary}
    \label{alg:boundary}
    \begin{algorithmic}[1]
\Require Octree $\mathcal{O}$, Element marker $\mathcal{M}$
\Ensure Surrogate boundary $\Tilde{\Gamma}$, Marker $\mathcal{M}$
\item[]

\For{element $\in \mathcal{O}$}
\State FaceBits $\leftarrow$ [false]
\For{face $\in$ element} \Comment{Loop over the faces of the element}
\If{$\mathcal{M}$[element] == \Intercepted{} and \texttt{allof}(nodes $\in$ face == \Exterior{})} \Comment{Condition for \Intercepted element}
\State FaceBits[face] $\leftarrow$ True
\ElsIf{$\mathcal{M}$[element] == \Neighbors{}} \Comment{Condition for \Neighbors elements}
\State BoundaryFace $\leftarrow$ True
\For{nodes $\in$ face}
\If{(nodes == \Exterior{}) or (nodes == \FalseIntercepted \textsc{Node})} \Comment{Each node should be either \Exterior{} or \FalseIntercepted }
\State continue
\Else 
\State BoundaryFace $\leftarrow$ False 
\State break
\EndIf
\EndFor
\EndIf
\EndFor
\If{cycle(faceBits)} \Comment{If opposite face of the same element forms the part of $\Tilde{\Gamma}$}
\State $\mathcal{M}$[element] $ \leftarrow$ \FalseIntercepted
\Else
\For{face $\in$ element} \Comment{Add faces to the surrogate domain $\Tilde{\Gamma}$}
\If{faceBits[face] $\leftarrow$ true}
\State $\Tilde{\Gamma} \leftarrow$ $\Tilde{\Gamma}$.\texttt{push\_back}(face)
\EndIf
\EndFor
\EndIf
\EndFor
\item[]
\Return ($\Tilde{\Gamma}$,$\mathcal{M}$)
\end{algorithmic}
\end{algorithm}

In some cases, we observe a cycle being formed where the opposite faces ($X_i^{-},X_i^{+}, i = 1 \cdots $ dim) of a given element are both chosen to be part of the surrogate boundary. This violates the second condition of the surrogate boundary described in \secref{sec: surrogate_math}. We mark such elements as \FalseIntercepted{} to resolve this issue. This ensures that only one of the sides $X_i^{-}$ or $X_i^{+}$ is selected to be part of the surrogate boundary depending on the side of \Neighbors{} (\Algref{alg:boundary}). \figref{fig:DifferenceLambda} shows the variation of the surrogate domain and surrogate boundary for different choices of $\lambda$. We note that all the steps in \Algref{alg:surrogateBoundary} can be done efficiently in $\mathcal{O}(N)$ steps and require a small number of passes over the elements of the octree.

\subsubsection{SBM computation}

Once we set up the surrogate domain ($\tO$) and resultant surrogate boundary ($\tG$) along with the associated markers, we proceed with the steps for the deployment of the SBM. \Algref{alg:sbm} briefs the major step for the SBM computation. We note that the elements marked as \FalseIntercepted{} (and \Exterior{}, if not using incomplete octree) are skipped over, whereas the volume integration is performed on other elements. Each face of the \Intercepted{} and \Neighbors{} element is checked to see if it belongs to the surrogate boundary $\tG$. If a given face belongs to $\tG$, the required surface integration (as given in~\eqnref{eq:SBM}) computation is performed and assembled to the global matrix and vector. Once the global matrix $\mathcal{A}$ and global vector $b$ are assembled, we solve the system of equations to obtain the solution $u^h$.

\begin{algorithm}[h!]
  \footnotesize
    \caption{\textsc{SBM:} Shifted Boundary method}
    \label{alg:sbm}
    \begin{algorithmic}[1]
\Require Octree $\mathcal{O}$, Element marker $\mathcal{M}$,Surrogate boundary $\Tilde{\Gamma}$, PDE ($\mathcal{L}(u) = f$)
\Ensure Solution of PDE: $u^h$
\item[]
\For{element $\in \mathcal{O}$}
\If{element == \FalseIntercepted} 
\State continue
\Else
\State Assemble matrix ($\mathcal{A}$) and right hand side vector ($b$)
\If {(element == \Intercepted) or (element == \Neighbors) }
\For{face $\in$ element}
\If{face $\in \Tilde{\Gamma}$}
\State Assemble surface contribution to matrix ($\mathcal{A}$) and right hand side vector ($b$)
\EndIf
\EndFor
\EndIf
\EndIf
\EndFor
\State Solve $\mathcal{A}u^h = b$ \Comment{Solve system of linear equation}
\item[]
\Return{$u^h$}
\end{algorithmic}
\end{algorithm}

\subsection{Distance function calculation}

Note that the SBM computations require evaluating the distance, $\bs{d}$, of Gauss points on the surrogate boundary to the closed point on the true boundary. This section focuses on calculating distance functions for intricate three-dimensional geometries. We consider the geometries to be represented in STL files, which are, in turn, represented by sets of triangles. SBM computation requires computing the distance function by finding the normal distance from the Gauss point to the nearest triangle. We store information about Gauss points and their corresponding distance functions to avoid repeated distance function calculations. This is particularly important for time-dependent problems on the static mesh as it prevents repetitive computation at the cost of extra memory.

\input{Tikz/DistanceCalc}
\begin{algorithm*}[t!]
\caption{\textsc{Overview}: Calculation of distance functions}
\footnotesize
\begin{algorithmic}[1]
\Require Gauss point positions on the surrogate boundary (GaussPoints), geometry, and mapping between Gauss points and distance functions
\Ensure Distance function corresponding to each Gauss point
\item[]
\For{$P \in $ GaussPoints}
\If{$P$ is already in the mapping}
\State Use the precomputed distance function for $P$
\Else
\State Find the distance function from $P$ to the nearest triangle using \textsc{NormalDistCalc} 
\Comment{~\Algref{Alg: NormalDistCalc}}
\State Save the mapping between $P$ and its distance function
\EndIf
\EndFor
\State \Return distance functions for all Gauss points
\end{algorithmic}
\label{Alg: Overview}
\end{algorithm*}

\begin{algorithm*}[t!]
\caption{\textsc{NormalDistCalc}: Calculation of normal distance to nearest triangle}
\label{Alg: NormalDistCalc}
\footnotesize
\begin{algorithmic}[1]
\Require k-d tree built from triangle centroids, Gauss point position ($P$), and geometry information
\Ensure Distance function corresponding to the Gauss point
\item[]
\State Use the k-d tree to find the ID of the nearest triangle to the Gauss point $P$
\State Calculate the normal vector from $P$ to the nearest triangle by projecting the vector $\mathbf{PA}$ onto the triangle's normal vector $\mathbf{n}$: $(\frac{\mathbf{PA} \cdot \mathbf{n}}{ \lvert \mathbf{n} \rvert ^2 }) \cdot \mathbf{n}$, where $A$ is one vertex position of the triangle 
\If{the projection point is inside the 3D triangle} \Comment{~\Algref{Alg: CheckInside3DTriangle}}
\State Set the distance function as the normal vector to the nearest triangle
\Else
\State Calculate the shortest distance from $P$ to the nearest triangle edge as the distance function
\Comment{~\Algref{Alg: ShortestDist2TriEdge}}
\EndIf
\State \Return distance function
\end{algorithmic}
\end{algorithm*}

\begin{algorithm*}[t!]
\caption{\textsc{CheckInside3DTriangle}: Check if point is inside 3D triangle}
\label{Alg: CheckInside3DTriangle}
\footnotesize
\begin{algorithmic}[1]
\Require Projection point position ($P_{project}$) and vertex positions of triangle ($A$, $B$, $C$)
\Ensure Whether the projection point is inside the 3D triangle
\item[]
\State Calculate cross products $\mathbf{u} = \mathbf{AB} \times \mathbf{AP_{project}},\mathbf{v} = \mathbf{BC} \times \mathbf{BP_{project}}, \mathbf{w} = \mathbf{CA} \times \mathbf{CP_{project}}$
\If{$\mathbf{u} \cdot \mathbf{v}$ < 0 or $\mathbf{u} \cdot \mathbf{w}$ < 0} 
\State \Return false
\Else
\State \Return true
\EndIf
\end{algorithmic}
\end{algorithm*}

\begin{algorithm*}[t!]
\caption{\textsc{ShortestDist2TriEdge}: Calculate the shortest distance to a triangle edge}
\label{Alg: ShortestDist2TriEdge}
\footnotesize
\begin{algorithmic}[1]
\Require Gauss point position $P$, and vertices of a triangle $A$, $B$, and $C$
\Ensure Shortest distance vector from $P$ to a triangle edge
\item[]
\State Compute the projection of $P$ onto each triangle edge, i.e., $\mathbf{PP_{AB}} = \frac{\mathbf{AP} \cdot \mathbf{AB}}{{\lvert \mathbf{AB} \rvert}^2} \mathbf{AB} - \mathbf{AP}$, $\mathbf{PP_{BC}} = \frac{\mathbf{BP} \cdot \mathbf{BC}} {{\lvert \mathbf{BC} \rvert}^2} \mathbf{BC} - \mathbf{BP}$, $\mathbf{PP_{CA}} = \frac{\mathbf{CP} \cdot \mathbf{CA}}{{\lvert \mathbf{CA} \rvert}^2} \mathbf{CA} - \mathbf{CP}$, where $P_{AB}$, $P_{BC}$, $P_{CA}$ denote the closest points on edges $\overline{AB}$, $\overline{BC}$, $\overline{CA}$ to $P$
\State Compute the minimum distance from $P$ to each edge, i.e., $\lvert \mathbf{PP_{AB}} \rvert$, $\lvert \mathbf{PP_{BC}} \rvert$, $\lvert \mathbf{PP_{CA}} \rvert$ and find the closest point ($P_{closest}$) to $P$ within $P_{AB},P_{BC},P_{CA}$
\If{the $P_{closest}$ is inside the 3D triangle} \Comment{~\Algref{Alg: CheckInside3DTriangle}}
\State Set the distance function as the $\mathbf{PP_{closest}}$
\Else
\State Calculate the shortest distance function from $P$ to the nearest triangle vertex, $\mathbf{PP_{vertex}}$ 
\EndIf
\State \Return distance function
\end{algorithmic}
\end{algorithm*}

The procedure for calculating the distance function is outlined in \Algref{Alg: Overview} and \figref{fig:FlowChart}. We utilize the algorithm presented in \Algref{Alg: NormalDistCalc} to compute the normal distance between Gauss points and the nearest triangle, as depicted in \figref{fig:FallIn}. We leverage a k-d tree to efficiently search the closest triangle using \texttt{nanoflann} library~\citep{blanco2014nanoflann} to expedite this process. Specifically, the k-d tree is constructed from the triangle centroids, and during the Gauss point iteration, the k-d tree assists in identifying the nearest triangle and obtaining its ID for each Gauss point. Nonetheless, in certain instances, the projection points from Gauss points to triangles may fall outside the triangle, as illustrated in \figref{fig:FallOut}. To handle such situations, we have incorporated two additional procedures. Firstly, \Algref{Alg: CheckInside3DTriangle} verifies whether the projection points lie within the nearest triangle. Secondly, \Algref{Alg: ShortestDist2TriEdge} determines the shortest distance between the point and the edges of the triangle, as shown in \figref{fig:ShortEdge}. It is worth noting that the projection method we use to find the closest point on the triangle edges to the Gauss point may not always result in a point inside the triangle in three dimensions, as shown in \figref{fig:ShortEdgeOutside}. In such cases, we search for the closest vertex of the triangle and calculate the distance function based on it, as illustrated in \figref{fig:2Vertex} and \Algref{Alg: ShortestDist2TriEdge}.

\pagebreak

\section{Numerical results}
\label{sec: Results}
This section presents numerical results for the simulation of Poisson equation and the equations of linear elasticity, over domains of complex geometry. We note the following important points for the readers to interpret the results:
\begin{itemize}
    \item The integration for the weak form is performed using standard $(p + 1)^{dim}$ Gauss quadrature points in all the elements, where $p$ is the order of the polynomial finite element interpolation basis. We use the linear basis function ($p = 1$) in all reported results.
    \item We report accuracy by comparing the numerical solution against analytical solutions. This post-processing operation to compute the $L^2$-error is performed on each element using five Gauss points per dimension.
    The reported $L^2$-error is computed as 
    $$    ||e||_{L^2(\Om)}
    =
    ||u^h - u_{exact}||_{L^2(\Om)} = \sqrt{\int_{\Om} (u^h - u_{exact})^2 d\Om} \; .
    $$
    Note that the error is computed on the true domain $\Om$. In particular, the SBM solution is smoothly extended over elements that intersect $\Om$ but are not part of the SBM active domain $\tO$. 
    \item In order to compare errors from different simulations, for which the sizes of the domains $\Om$ may differ, we report the normalized error, defined as:    $$ L_{2N} (\Om) = \frac{||e||_{L^2(\Om)}}{\sqrt{\int_{\Om}d\Om}} \; . $$
    \item We define the improvement factor $I_{2\lambda}$ (\eqnref{eq:l2lambda}) as the metric for comparison between different surrogates with respect to the case in which $\lambda = 1$. We recall that $\lambda=1$ corresponds to the case in which $\tO \supset \Om$, that is $\tO$ circumscribes $\Om$.    
    When comparing simulations performed with different values of $\lambda$, a lower value of $I_{2\lambda}$, corresponds to better solutions:
    $$
    I_{2\lambda} = \frac{L_{2N}(\Om;\lambda)}{L_{2N}(\Om;\lambda = 1)}
        \label{eq:l2lambda}
        \; .
    $$
    \item We recall the definition of $\lambda, (0 \leq \lambda \leq 1)$ as the elemental volume fraction of the domain $\Om$ . More specifically, $\lambda$ is a threshold used to check whether an \Intercepted{} element needs also to be classified as  \FalseIntercepted{} element.
\end{itemize}

\subsection{Identifying the optimal surrogate boundary}
Recall that our goal is to construct the \emph{optimal} surrogate boundary that minimizes the distance $\bs{d}$, the distance between $\Gamma$ and $\tG$. Our approach yields that $\lambda = 0.5$ minimizes $\bs{d}$. We illustrate this result using a canonical test example. Consider a rotated rectangle inclined at a 15$^{\circ}$ angle with respect to the background octree mesh. We compute the distance (\eqnref{eq: optim}, the RMS distance) between the surrogate boundary and true boundary for a range of surrogate boundary choices by varying $\lambda \in [0, 1]$.

\begin{figure}[b!]
    \centering
    \includegraphics[width=0.9\linewidth,trim=0 0 0 0,clip]      {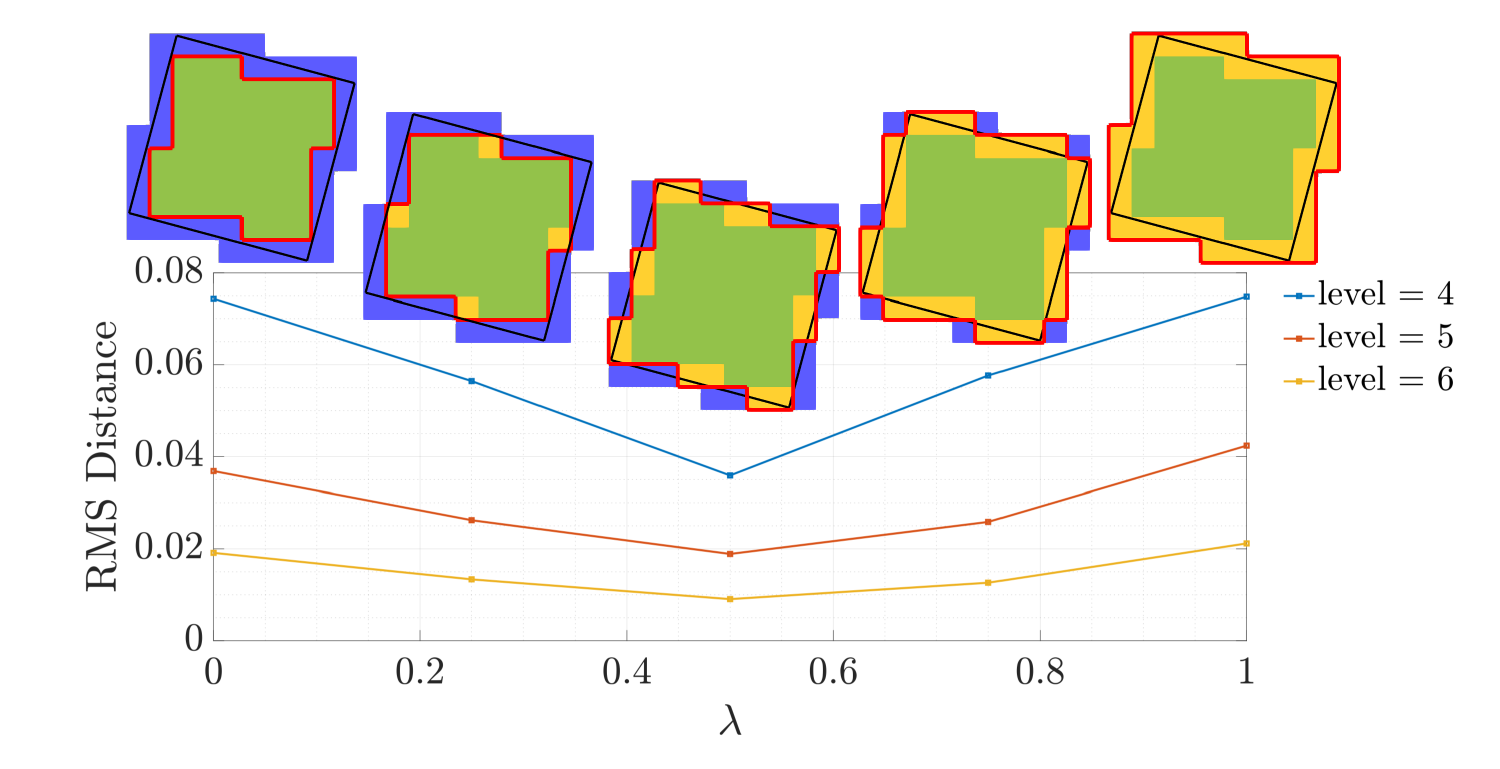}
    \caption{Constructing the surrogate boundary for a rotated square on a Cartesian mesh. The figure shows the RMS distance with varying $\lambda$ values for different octree levels. In each inset, the solid red line represents the surrogate surface, while the solid black line represents the true boundary. \FalseIntercepted{} elements are marked by \textcolor{blue}{$\blacksquare$}, whereas \TrueIntercepted{} elements are marked by \textcolor{cpu2}{$\blacksquare$} and \Interior{} elements by \textcolor{cpu1}{$\blacksquare$}.
    }
    \label{fig:IdentifyingOptSug}
\end{figure}

\begin{figure}[t!]
    \centering
    \includegraphics[width=0.8\linewidth,trim=0 0 0 0,clip]  {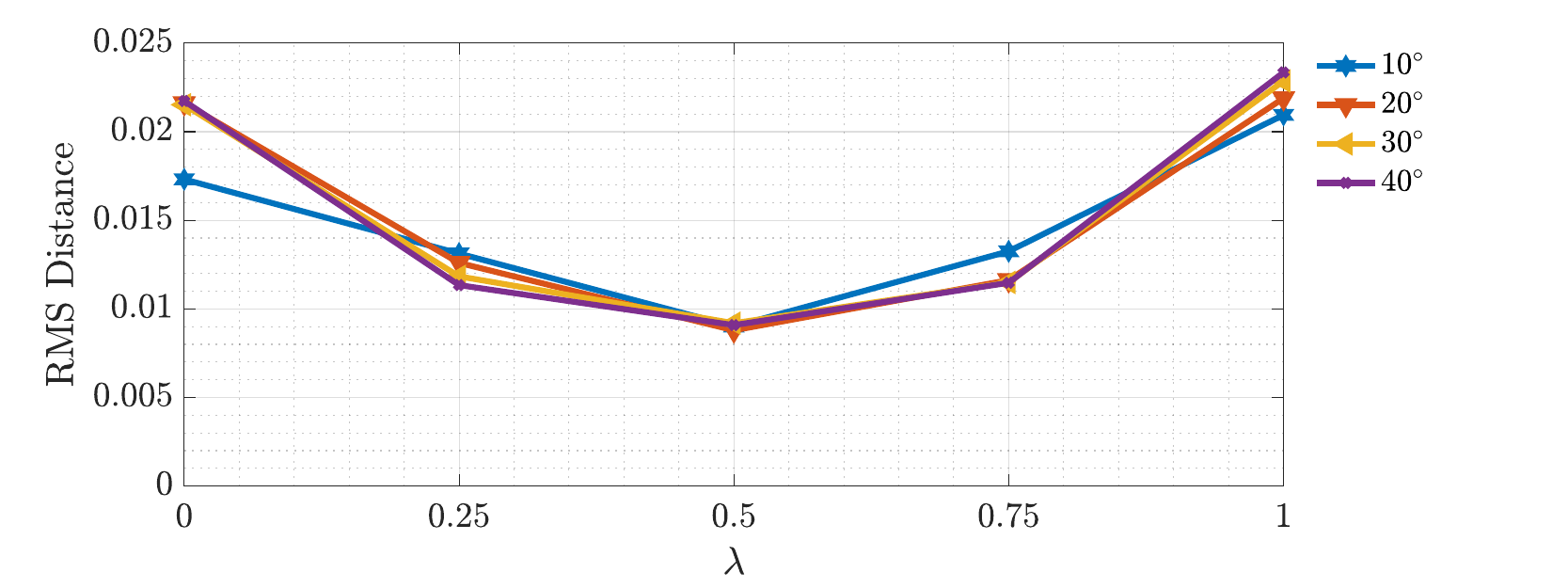}
    \caption{RMS of the distance between the true and surrogate boundaries, for a square rotated by various angles ($10^\circ$, $20^\circ$, $30^\circ$, $40^\circ$), and grid refinement level equal to 6.}
    \label{fig:IdentifyingOptSug_MultipleAngles}
\end{figure}

\figref{fig:IdentifyingOptSug} compares the value of RMS distance for different $\lambda$, for different mesh resolutions (identified by the level of refinement). Observe that $\lambda = 0.5$ yields the minimal RMS distance, irrespective of the level of the mesh resolution. In addition to conducting tests at a 15$^{\circ}$ angle, we also carried out experiments at 10$^{\circ}$, 20$^{\circ}$, 30$^{\circ}$, and 40$^{\circ}$ angles, paired with varying mesh refinement levels. The results consistently pointed towards $\lambda = 0.5$ as the value that assures the least RMS distance. \figref{fig:IdentifyingOptSug_MultipleAngles} further demonstrates these conclusions for various geometries defined by rotating a square at various angles. Thus, choosing $\lambda = 0.5$ provides a simple and effective strategy for constructing the surrogate boundary. In the subsequent sections, we analyze improvement in solution accuracy resulting from constructing an \emph{optimal} surrogate boundary ($\lambda = 0.5$) as compared to a non-optimal surrogate (usually $\lambda = 0, \lambda = 1$).

\subsection{Solving Poisson's equation on disk}
\label{subsec:PoissonDisk}
We next utilize the optimal surrogate construction approach to solve Poisson's equation when $\Om$ has a simple shape: a circular disk. This is an exact geometry, and several ways exist to construct an axis-aligned surrogate boundary for a circle. Consider a disk with a radius $R = 0.5$, centered at ($x_0 = 0.5, y_0 = 0.5$). We solve Poisson's equation on this geometry with a Dirichlet boundary of $u_0 = 0.01$ on the boundary and forcing term $f = 1$. We choose the penalty parameter $\alpha$ to be 400. This problem has an analytical solution given by: $u(r) = 0.25(R^2-r^2)+u_0$, where $r$ is the distance from the center ($r = \sqrt{(x - x_0)^2 + (y - y_0)^2}$).

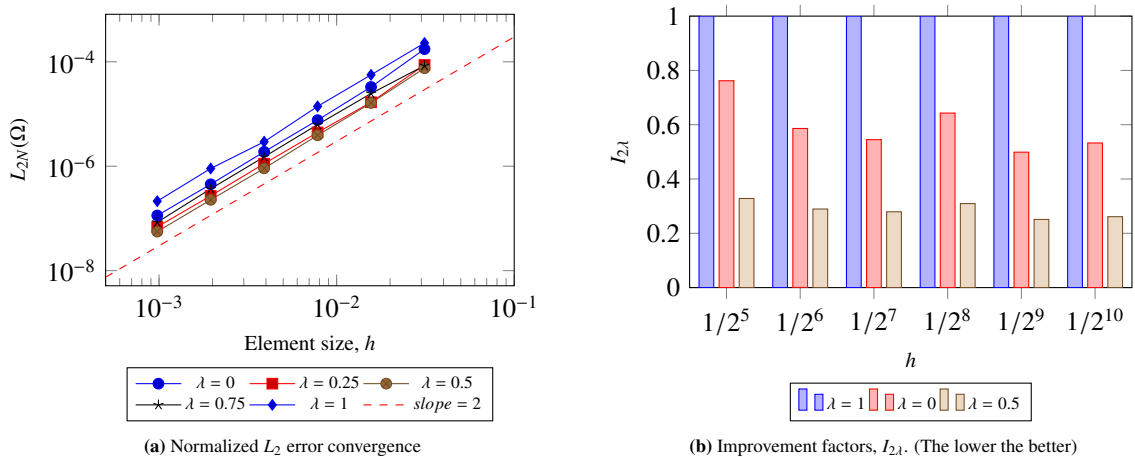
\begin{figure}[b!]
    \centering
    \begin{subfigure}{0.45\textwidth}
    \begin{tikzpicture}
    \begin{loglogaxis}[
        width=0.94\linewidth, 
        height = 0.7\linewidth,
        scaled y ticks=true,
        xlabel={\footnotesize Element size, $h$},
        ylabel={\footnotesize $L_{2N} (\Omega)$},
        legend entries={$\lambda = 0$, $\lambda = 0.25$, $\lambda = 0.5$, $\lambda = 0.75$, $\lambda = 1$, $slope = 2$},
        legend style={at={(0.5,-0.3)},anchor=north, nodes={scale=0.65, transform shape}}, 
        legend columns=3,
        xmin=5e-4,
        xmax=0.1	
        ]
        \addplot table [x={h},y={NormalizedError_0p0_adaptive},col sep=comma] {Figures/L2error/L2error2DcircleOpt.txt};
        \addplot table [x={h},y={NormalizedError_0p25_adaptive},col sep=comma] {Figures/L2error/L2error2DcircleOpt.txt};
        \addplot table [x={h},y={NormalizedError_0p5_adaptive},col sep=comma] {Figures/L2error/L2error2DcircleOpt.txt};
        \addplot table [x={h},y={NormalizedError_0p75_adaptive},col sep=comma] {Figures/L2error/L2error2DcircleOpt.txt};
        \addplot table [x={h},y={NormalizedError_1p0_adaptive},col sep=comma] {Figures/L2error/L2error2DcircleOpt.txt};
        \addplot +[mark=none, red, dashed] [domain=1e-4:1]{0.03*x^2};
    \end{loglogaxis}
    \end{tikzpicture}
    \caption{Normalized $L_2$ error convergence}
    \label{fig:Poisson_disk_NormL2}
    \end{subfigure}
    \hspace{0.02\linewidth}
    \begin{subfigure}{0.45\textwidth}
    \begin{tikzpicture}
    \begin{axis}[
        width = \linewidth,
        height = 0.7\linewidth,
        ybar,
        ymax = 1.0,
        ymin= 0.0,
        xtick= {0,1,2,3,4,5},
        xticklabels={$1/2^5$,$1/2^6$,$1/2^7$,$1/2^8$,$1/2^9$,$1/2^{10}$},
        bar width=0.2,
        xlabel={\footnotesize $h$},
        ylabel = {\footnotesize $I_{2\lambda}$},
        legend style={at={(0.5,-0.35)},anchor=north, nodes={scale=0.65, transform shape}}, 
        legend columns=3,
    ]
    \addplot
    +[]
	table[x expr=\coordindex,y expr={\thisrow{NormalizedError_1p0_adaptive}/\thisrow{NormalizedError_1p0_adaptive}},col sep=comma]{Figures/L2error/L2error2DcircleOpt.txt};
 
    \addplot
    +[]
	table[x expr=\coordindex,y expr={\thisrow{NormalizedError_0p0_adaptive}/\thisrow{NormalizedError_1p0_adaptive}},col sep=comma]{Figures/L2error/L2error2DcircleOpt.txt};
 
	\addplot
    +[]
	table[x expr=\coordindex,y expr={\thisrow{NormalizedError_0p5_adaptive}/\thisrow{NormalizedError_1p0_adaptive}},col sep=comma]{Figures/L2error/L2error2DcircleOpt.txt};
	\legend{$\lambda=1$,$\lambda=0$,$\lambda=0.5$}
    \end{axis}
    \end{tikzpicture}
    \caption{Improvement factors, $I_{2\lambda}$. (The lower the better)
    }
    \label{fig:Poisson_disk_improvement}
\end{subfigure}
\caption{Left: Mesh convergence plot solving the Poisson's equation on a disk. Each line represents a convergence plot with a specific surrogate boundary constructed using a $\lambda$. Notice that while all surrogate choices exhibit the expected slope, the $\lambda = 0.5$ surrogate produces the lowest error for any mesh size. Right: The improvement plot shows that the optimal surrogate produces more accurate solutions.}
\label{fig:Poisson_disk}
\end{figure}

\figref{fig:Poisson_disk} shows the mesh convergence  plot for different choices of $\lambda$. We observe second-order convergence in normalized error  $L_{2N}$ (\figref{fig:Poisson_disk_NormL2}), irrespective of the choice of $\lambda$. However, the optimal choice of $\lambda$, and thereby of surrogate boundary, can significantly reduce the magnitude of the error. We observe that choosing $\lambda = 0.5$ yields minimal error that is almost an order of magnitude lower when compared to the previously reported surrogate boundary choices of $\lambda=0$~\citep{main2018shifted} or $\lambda=1$~\citep{saurabh2021scalable}. \figref{fig:Poisson_disk_improvement} compares this improvement factor for the cases of $\lambda=0.0$ and $\lambda=0.5$ with respect to case of $\lambda=1$ (note that $I_{2\lambda}=1$ for $\lambda=1$). The choice $\lambda=0.5$ produces errors at least three times lower than the case $\lambda = 0$ or $\lambda = 1$ across all mesh resolutions. This improvement is significant, considering that the computational effort in identifying the optimal surrogate is minimal.

\subsection{Complex geometries}

We next showcase the ability of the SBM to accurately solve the Poisson's equation over three-dimensional objects. We consider three standard benchmarks exhibiting complex geometries and sharp corners: the  Stanford Bunny, the Moai, and the Armadillo. For all three-dimensional cases, the penalty parameter $\alpha$ is set to 400. We use the method of manufactured solutions to construct an analytical solution against which to compare the SBM results. This allows rigorous comparison across different values of $\lambda$. The analytical solution for each of the cases is given as:
\begin{equation}  \label{eq:3DManufacturedSolutions}
u(x,y,z) =\left\{
\begin{aligned}
&\cos(\pi x) \, y \,  \sin(\pi z),&  \mathrm{\;Stanford \; bunny;} \\
&(1-x)\, (1-y)\, \cos(3 \pi z),   & \mathrm{\;Moai;} \\
&\cos(\pi x) \, (1-y) \, \sin(\pi z), & \mathrm{\;Armadillo.} 
\end{aligned}
\right.
\end{equation}

\begin{figure}[t!]
\centering
\begin{tabular}[t]{|c|c|c|}
\hline
\multicolumn{3}{|c|}{\small Stanford Bunny}
\\
\hline
    \includegraphics[width=0.25\linewidth,trim=650 0 650 0,clip, margin=0pt 0pt 0pt 1mm]{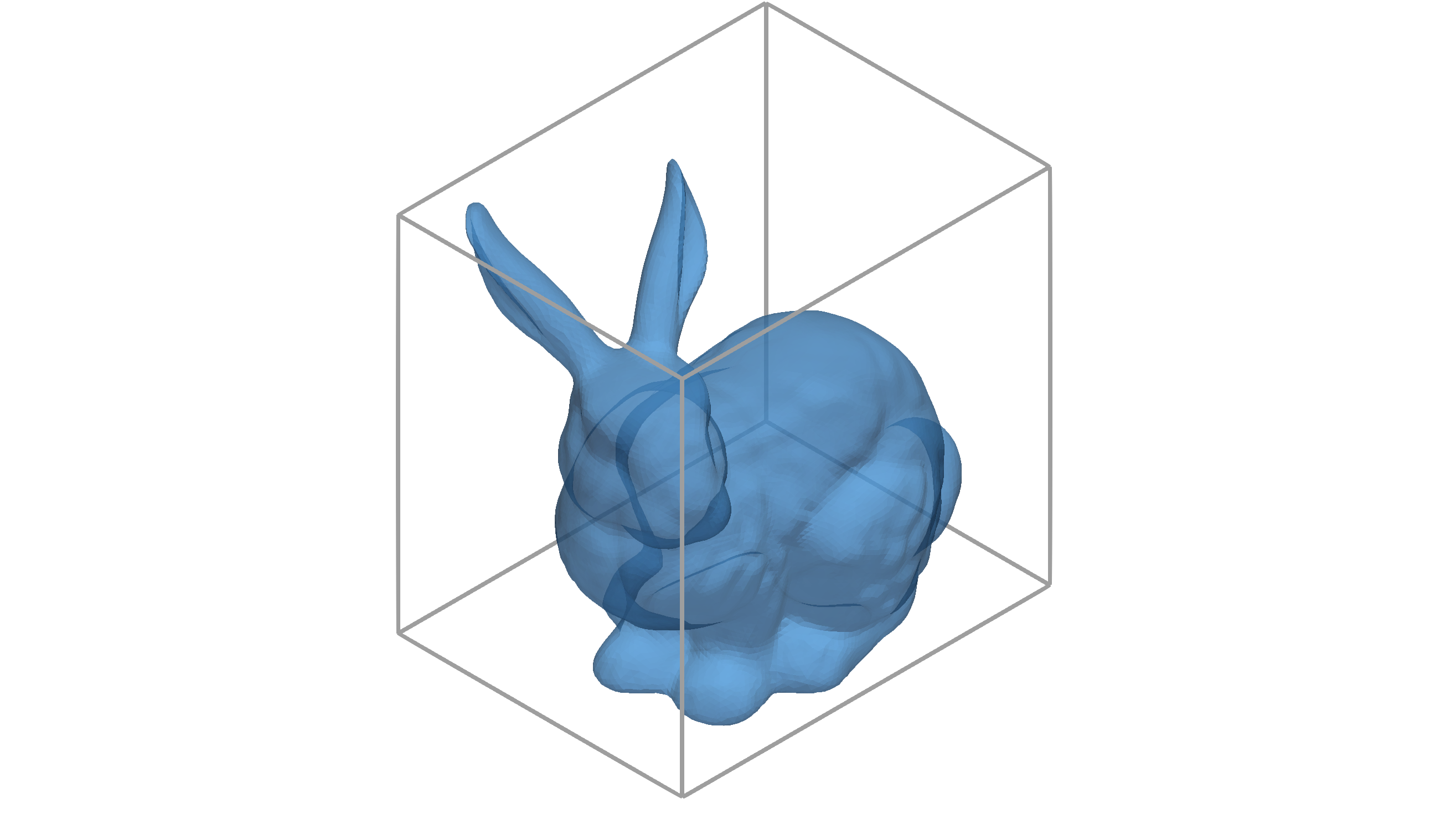}
    &
    \includegraphics[width=.25\linewidth,trim=80 100 80 95,clip]{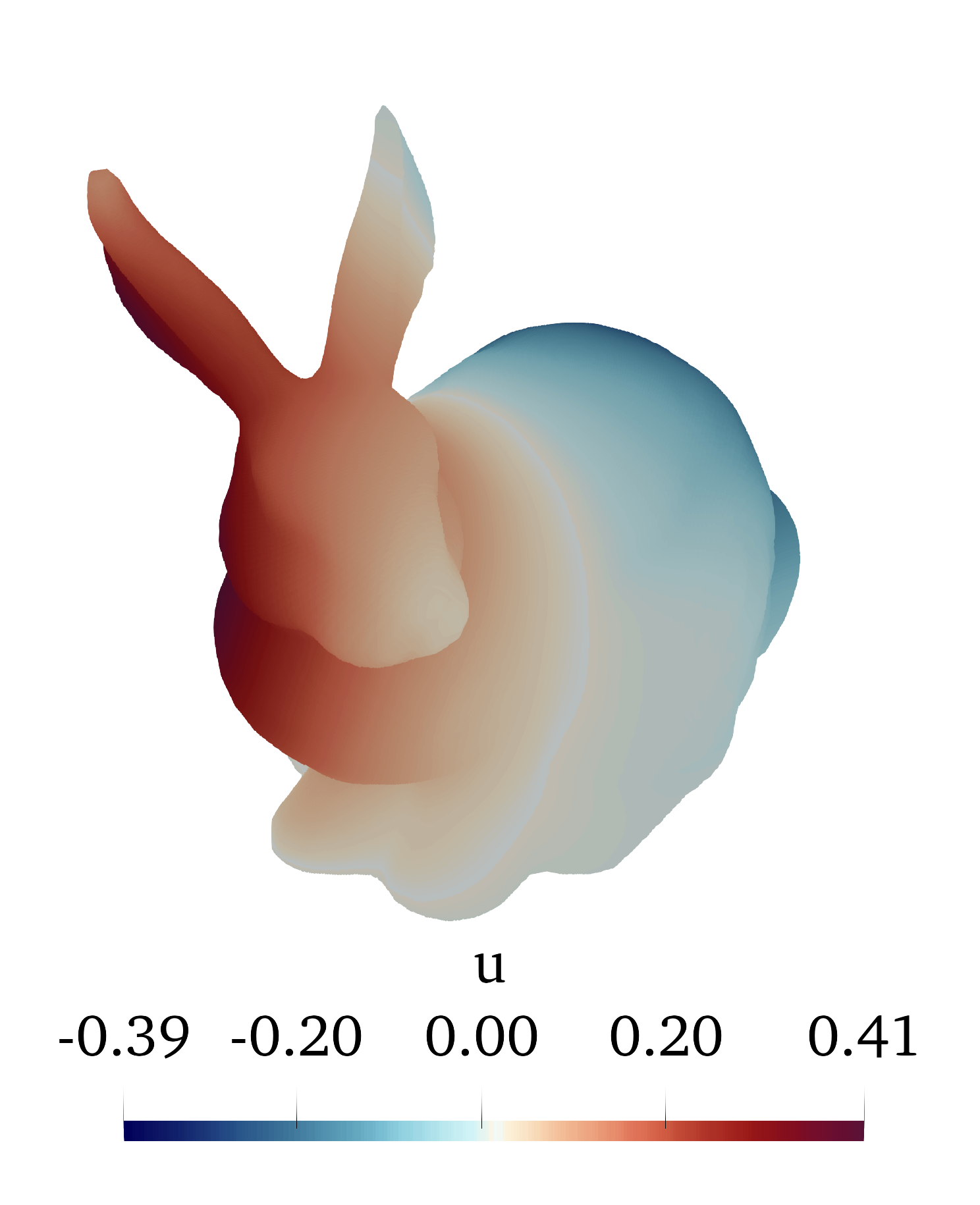}
    &
    \begin{tikzpicture}
    \begin{axis}[
        width = 0.35\linewidth,
        height = 0.25\linewidth,
        ybar,
        xtick= {0,1,2,3,4},
        xticklabels={$1/2^4$,$1/2^5$,$1/2^6$,$1/2^7$,$1/2^8$},
        bar width=0.2,
        xlabel = {$h$}, 
        ylabel = {$I_{2\lambda}$},
        legend style={at={(0.5,1.25)},anchor=north, nodes={scale=0.65, transform shape}}, legend columns=3,
    ]
    \addplot
     +[]
    	table[x expr=\coordindex,y expr={\thisrow{NormalizedError_1p0_adaptive}/\thisrow{NormalizedError_1p0_adaptive}},col sep=comma]{Figures/L2error/L2error3DbunnyOpt.txt};
    \addplot
     +[]
    	table[x expr=\coordindex,y expr={\thisrow{NormalizedError_0p0_adaptive}/\thisrow{NormalizedError_1p0_adaptive}},col sep=comma]{Figures/L2error/L2error3DbunnyOpt.txt};
 
    \addplot
    +[]
        table[x expr=\coordindex,y expr={\thisrow{NormalizedError_0p5_adaptive}/\thisrow{NormalizedError_1p0_adaptive}},col sep=comma]{Figures/L2error/L2error3DbunnyOpt.txt};
    \legend{$\lambda=1$,$\lambda=0$,$\lambda=0.5$}
    \end{axis}
\end{tikzpicture}
\\
\hline
\multicolumn{3}{|c|}{\small Moai}
\\
 \hline
    \includegraphics[width=.2\linewidth,trim=800 0 700 8,clip, margin=0pt 0pt 0pt 1mm]{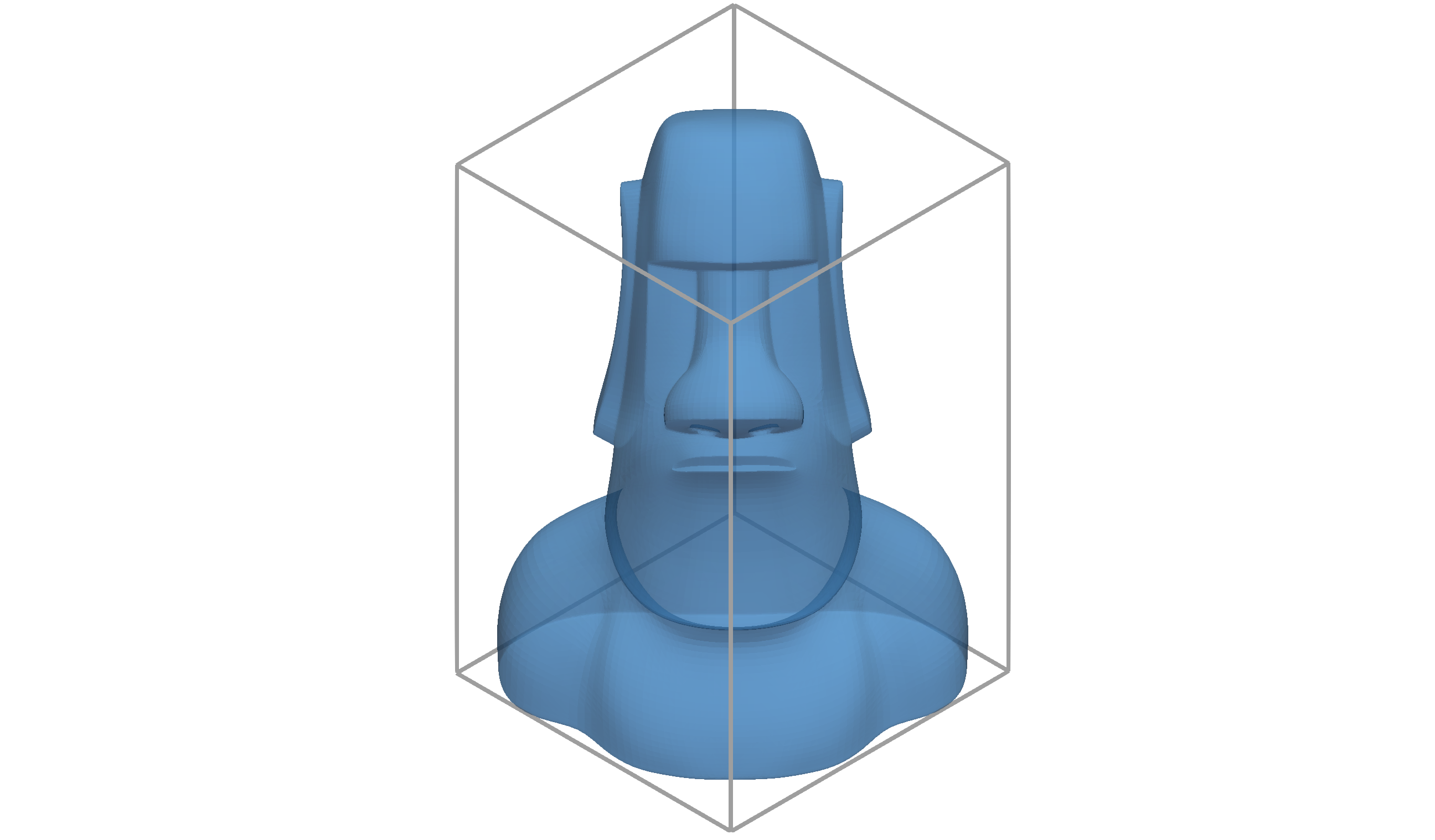}
    &
    \includegraphics[width=0.27\linewidth,trim=0 70 0 120,clip]{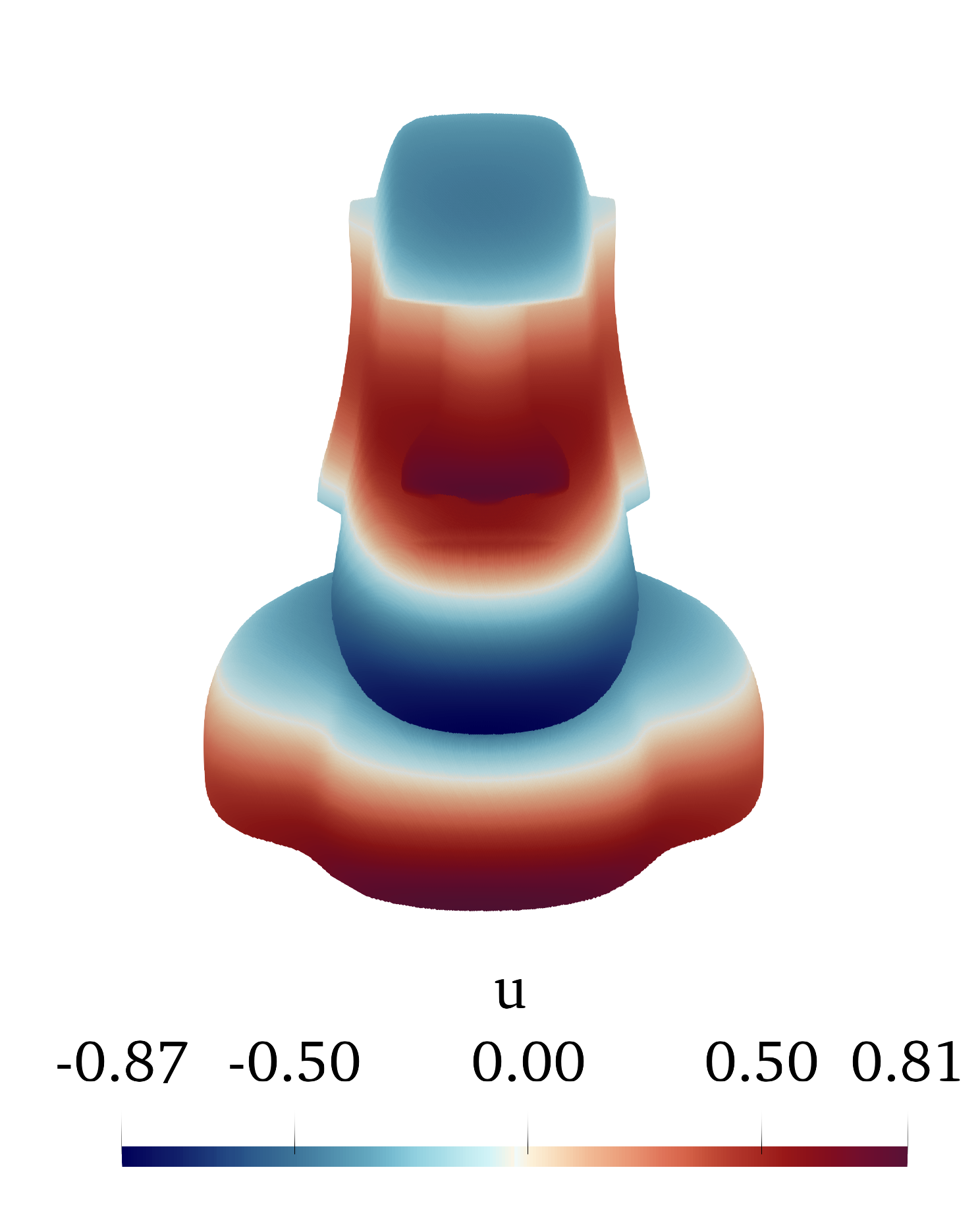}
    &
    \begin{tikzpicture}
    \begin{axis}[
        width = 0.35\linewidth,
        height = 0.25\linewidth,
        ybar,
        xtick= {0,1,2,3,4},
        xticklabels={$1/2^4$,$1/2^5$,$1/2^6$,$1/2^7$,$1/2^8$},
        bar width=0.2,
        xlabel = {$h$}, 
        ylabel = {$I_{2\lambda}$},
        legend style={at={(0.5,1.25)},anchor=north, nodes={scale=0.65, transform shape}}, legend columns=3,
    ]
    \addplot
     +[]
    	table[x expr=\coordindex,y expr={\thisrow{NormalizedError_1p0_adaptive}/\thisrow{NormalizedError_1p0_adaptive}},col sep=comma]{Figures/L2error/L2error3DmoaiOpt.txt};
     
    \addplot
     +[]
    	table[x expr=\coordindex,y expr={\thisrow{NormalizedError_0p0_adaptive}/\thisrow{NormalizedError_1p0_adaptive}},col sep=comma]{Figures/L2error/L2error3DmoaiOpt.txt};
     
    	\addplot
     +[]
    	table[x expr=\coordindex,y expr={\thisrow{NormalizedError_0p5_adaptive}/\thisrow{NormalizedError_1p0_adaptive}},col sep=comma]{Figures/L2error/L2error3DmoaiOpt.txt};
    \legend{$\lambda=1$,$\lambda=0$,$\lambda=0.5$}
    \end{axis}
\end{tikzpicture}

\\
\hline
\multicolumn{3}{|c|}{\small Armadillo}
\\
\hline
    \includegraphics[width=0.28\linewidth,trim=800 0 650 0,clip]{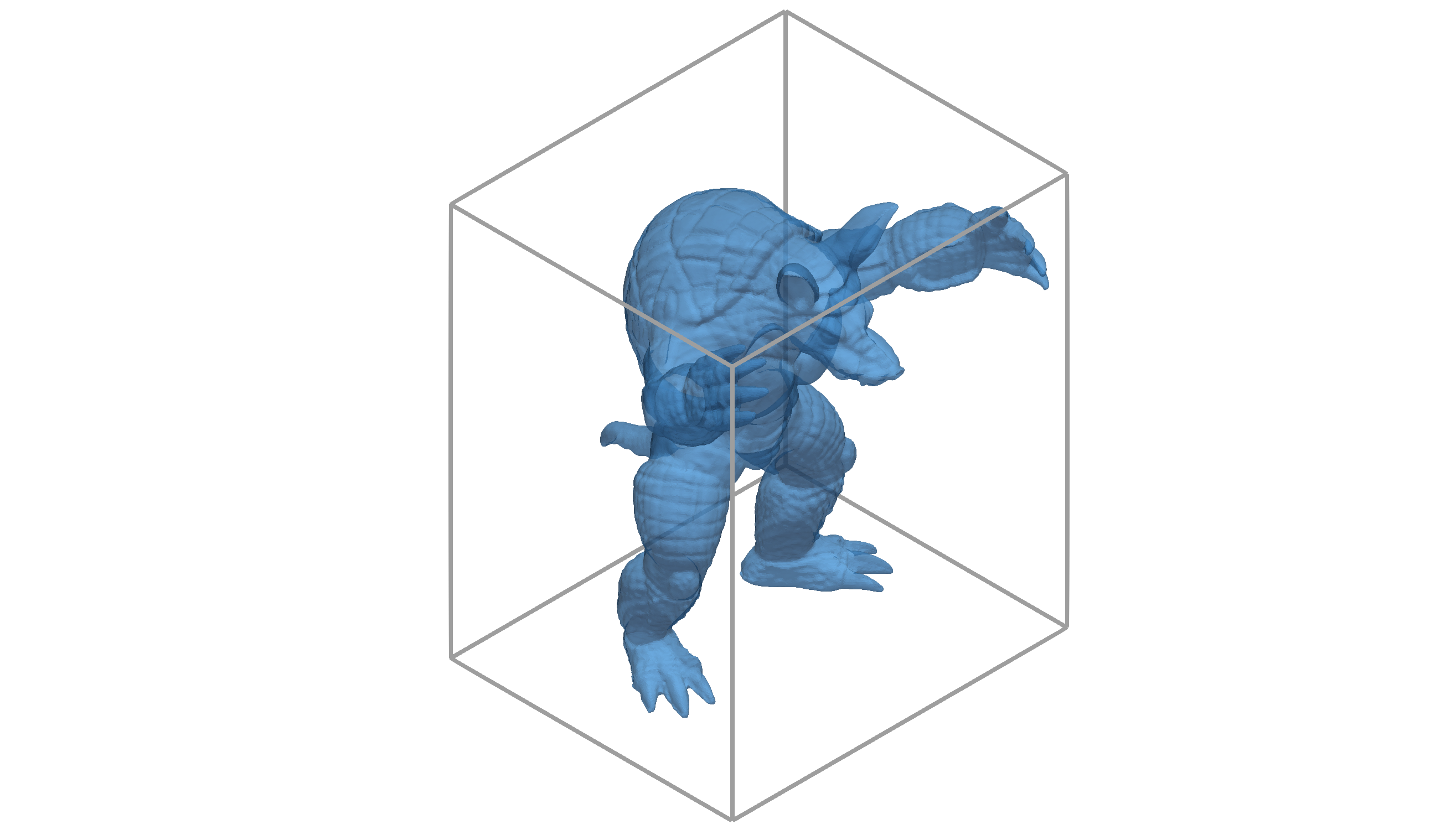}
    &
    \includegraphics[width=0.25\linewidth,trim=0 80 0 70,clip]{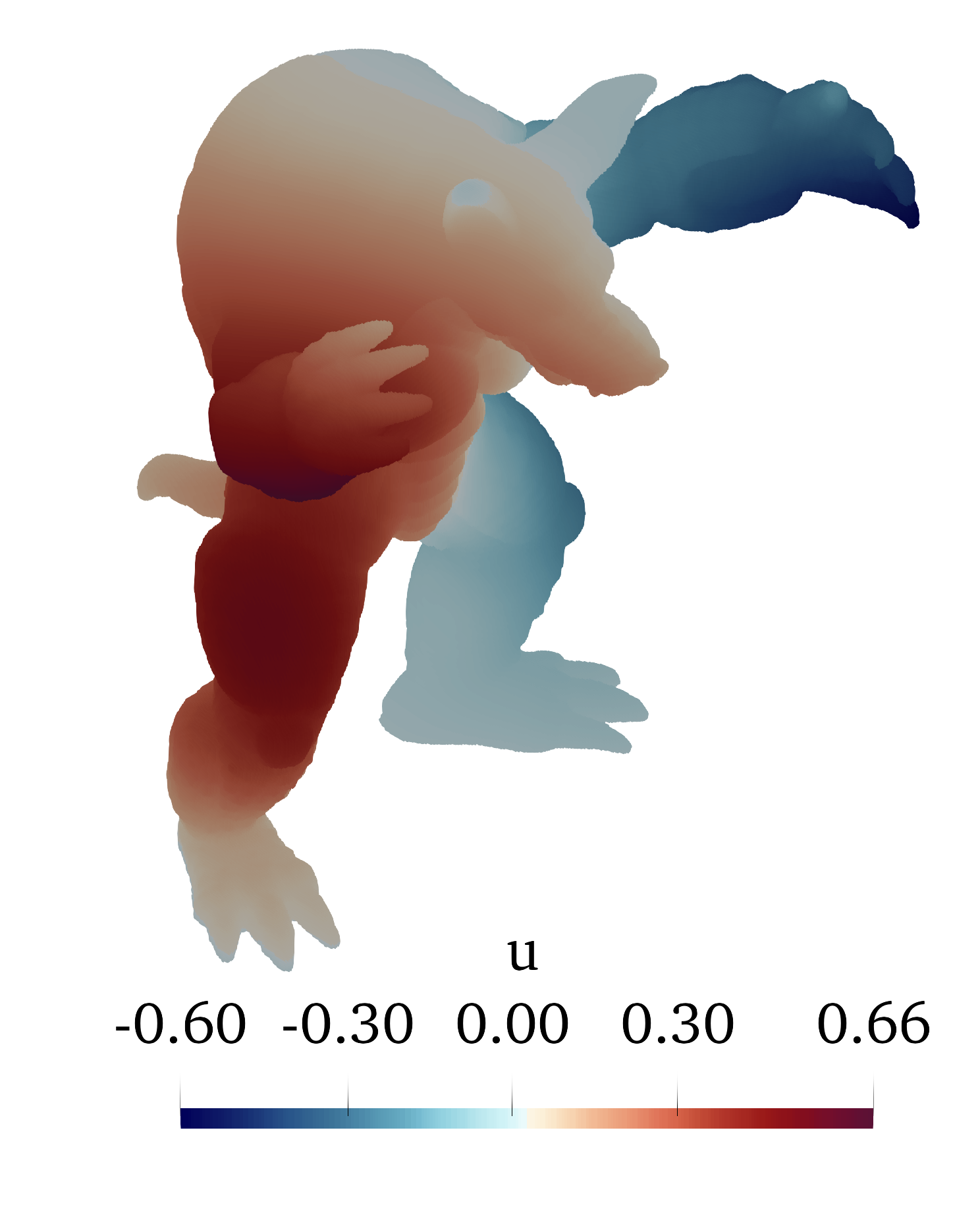}
    &
    \begin{tikzpicture}
    \begin{axis}[
        width = 0.35\linewidth,
        height = 0.25\linewidth,
        ybar,
        xtick= {0,1,2,3,4},
        xticklabels={$1/2^4$,$1/2^5$,$1/2^6$,$1/2^7$,$1/2^8$},
        bar width=0.2,
        xlabel = {$h$}, 
        ylabel = {$I_{2\lambda}$},
        legend style={at={(0.5,1.25)},anchor=north, nodes={scale=0.65, transform shape}}, legend columns=3,
    ]
    \addplot
     +[]
    	table[x expr=\coordindex,y expr={\thisrow{NormalizedError_1p0_adaptive}/\thisrow{NormalizedError_1p0_adaptive}},col sep=comma]{Figures/L2error/L2error3DarmOpt.txt};
    \addplot
     +[]
    	table[x expr=\coordindex,y expr={\thisrow{NormalizedError_0p0_adaptive}/\thisrow{NormalizedError_1p0_adaptive}},col sep=comma]{Figures/L2error/L2error3DarmOpt.txt};
	\addplot
    +[]
        table[x expr=\coordindex,y expr={\thisrow{NormalizedError_0p5_adaptive}/\thisrow{NormalizedError_1p0_adaptive}},col sep=comma]{Figures/L2error/L2error3DarmOpt.txt};
	\legend{$\lambda=1$,$\lambda=0$,$\lambda=0.5$}
    \end{axis}
\end{tikzpicture}
\\
\hline

\end{tabular}
\caption{Solving Poisson's equations on complex 3D domains using the optimized surrogate boundary.} 
\label{fig:3Dcase}
\end{figure}

For each object, we perform a mesh convergence analysis by solving the Poisson's equation using the optimal surrogate ($\lambda = 0.5$) and compare the accuracy against the choice of surrogates with $\lambda = 1$ and $\lambda = 0$. \figref{fig:3Dcase} compares $I_{2\lambda}$ for various values of $\lambda$ and different mesh sizes. We can see that the value of $\lambda = 0.5$ consistently outperforms $\lambda = 0$ and $\lambda = 1$ in terms of accuracy. This demonstrates both the importance of choosing an optimal $\lambda$ as well as the robustness of the proposed algorithm for solving PDEs on complex geometries.

\begin{figure}[t!]
    \centering
    \begin{subfigure}{0.3\textwidth}
        \centering
        \includegraphics[width=1.0\linewidth,trim=300 50 300 30,clip]{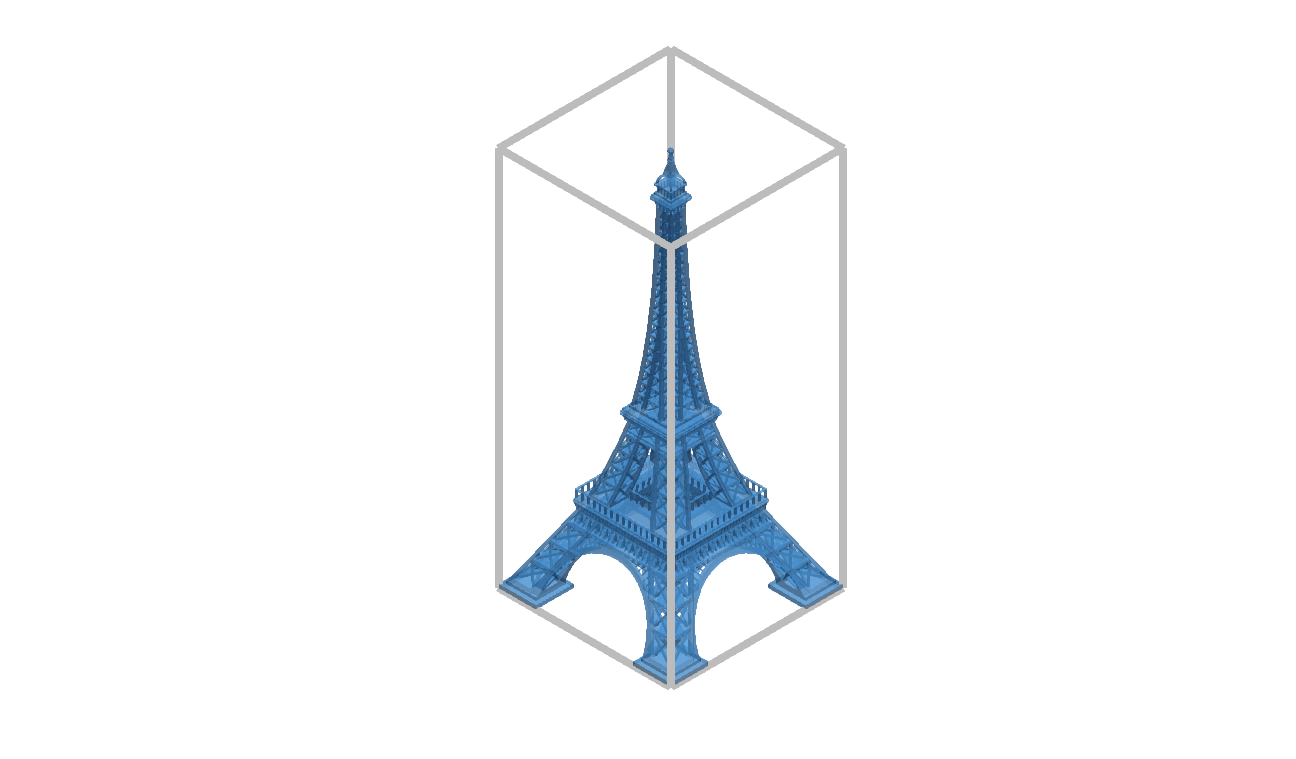}
        \caption{STL file}
        \label{fig:EiffelSTL}
    \end{subfigure}
    \hspace{0.01\linewidth}
    \begin{subfigure}{0.23\textwidth}
        \centering
        \includegraphics[width=1.0\linewidth,trim=0 0 0 0,clip]     {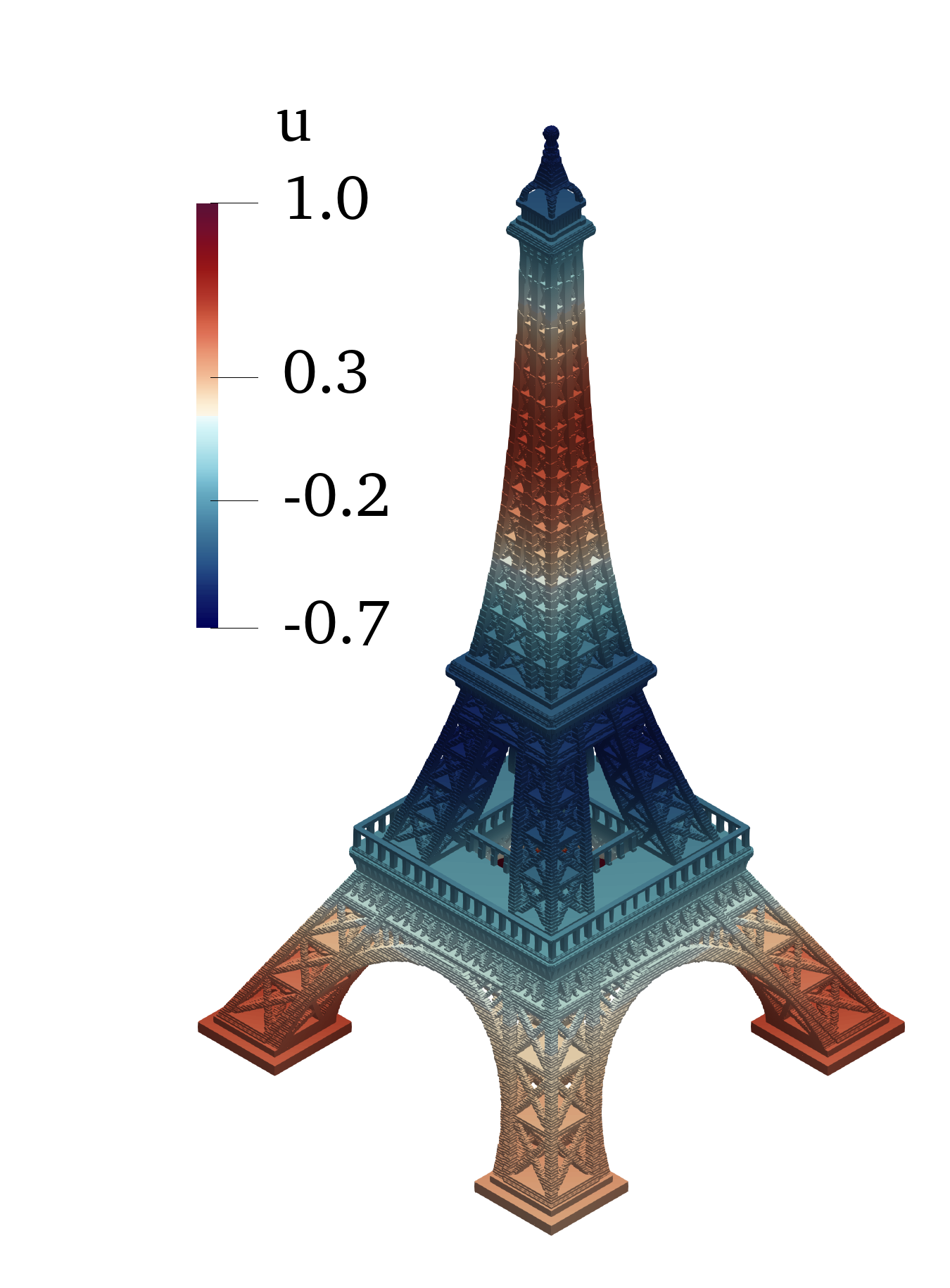}
        \caption{Solution contour}
        \label{fig:EiffelSolution}
    \end{subfigure}
    \hspace{0.01\linewidth}
    \begin{subfigure}{0.35\textwidth}
        \begin{tikzpicture}       
        \begin{axis}[
            width = \linewidth,
            height = 0.8\linewidth,
            ybar,
            xtick= {0,1,2,3,4,5},
            xticklabels={$1/2^{10}$},
            bar width=0.2,
            xlabel = {\footnotesize $h$},
            ylabel = {\footnotesize $I_{2\lambda}$},
            legend style={at={(0.5,1.25)},anchor=north, nodes={scale=0.65, transform shape}},
            legend columns=3,
            ]
    
    \addplot
     +[]
    	table[x expr=\coordindex,y expr={\thisrow{NormalizedError_1p0_adaptive}/\thisrow{NormalizedError_1p0_adaptive}},col sep=comma]{Figures/L2error/L2error3DeiffelOpt.txt};
     
    \addplot
     +[]
    	table[x expr=\coordindex,y expr={\thisrow{NormalizedError_0p0_adaptive}/\thisrow{NormalizedError_1p0_adaptive}},col sep=comma]{Figures/L2error/L2error3DeiffelOpt.txt};
     
    \addplot
     +[]
    	table[x expr=\coordindex,y expr={\thisrow{NormalizedError_0p5_adaptive}/\thisrow{NormalizedError_1p0_adaptive}},col sep=comma]{Figures/L2error/L2error3DeiffelOpt.txt};
    	\legend{$\lambda=1$,$\lambda=0$,$\lambda=0.5$}
       \end{axis}
    \end{tikzpicture}
    \caption{Error in the $L_2$-norm}
    \label{fig:EiffelL2}
\end{subfigure}
\caption{Solving Poisson's equation on the Eiffel tower model using SBM:(\textit{a}) The STL file of the Eiffel tower immersed into the background mesh and the bounding box. (\textit{b}) Solution contour of Eiffel Tower with mesh size equal to $1/2^{10}$. (\textit{c}) Improvement plot of $I_{2\lambda}$ and the results show that the optimal surrogate boundary can reduce the $L_{2N}$ error.}
\label{fig:EiffelTower}
\end{figure}

Next, to test the algorithm's robustness, we explore its performance on a geometry exhibiting an extremely complex topology. We consider a simulation on a three-dimensional model of the Eiffel Tower, in STL format. \figref{fig:EiffelSTL} shows the surface representation of the geometry, with a very large number of small holes and sharp corners. We choose a manufactured solution of the form:
\begin{equation}  \label{EiffelTower}
u(x,y,z) = (1-x)(1-y)\cos(3 \pi z) .
\end{equation}
\figref{fig:EiffelSolution} shows the resulting solution field. Finally, we notice the same trend in $L^2$-error, with surrogates defined by $\lambda = 0.5$ outperforming the other choices (\figref{fig:EiffelL2}). 

\input{Tikz/starLE.tex}

\subsection{Linear elasticity}

Our final example showcases the SBM approach for solving linear elasticity equations. As shown in \figref{fig:starLESoln}, we consider a star-shaped domain. The sharp corners and non-convex geometry makes this a challenging case. We set Young's modulus $E$ to be 1 and Poisson's ratio $\nu$ to be 0.3. The elastic tensor $\bs{C}$ is for plain stress. The penalty parameter $\gamma$ is equal to 400. We consider a manufactured solution of the form:
\begin{equation}  \label{eq:8}
\left\{
\begin{aligned}
u_{x}(x,y) =\frac{\sin(\pi x)\cos(\pi y)}{10} \\
u_{y}(x,y) =\frac{\cos(\pi x)\sin(\pi y)}{10}
\end{aligned}_{\textstyle \raisebox{2pt}{.}}\quad
\right.
\end{equation}
\figref{fig:starLESoln} shows the displacement solution contour, while \figref{fig:starLEError} shows the convergence of $L_{2N}$ error. Similar to Poisson's case, we observe a second-order convergence in $L^2$-error. We see a significant improvement in $L^2$-error by choosing a value of $\lambda = 0.5$ when compared to $\lambda = 0$ or $\lambda = 1$ both for the x-displacement, $u_x$ (\figref{fig:starLEux}) and y-displacement, $u_y$(\figref{fig:starLEuy}). This improved performance is seen across all mesh refinements.

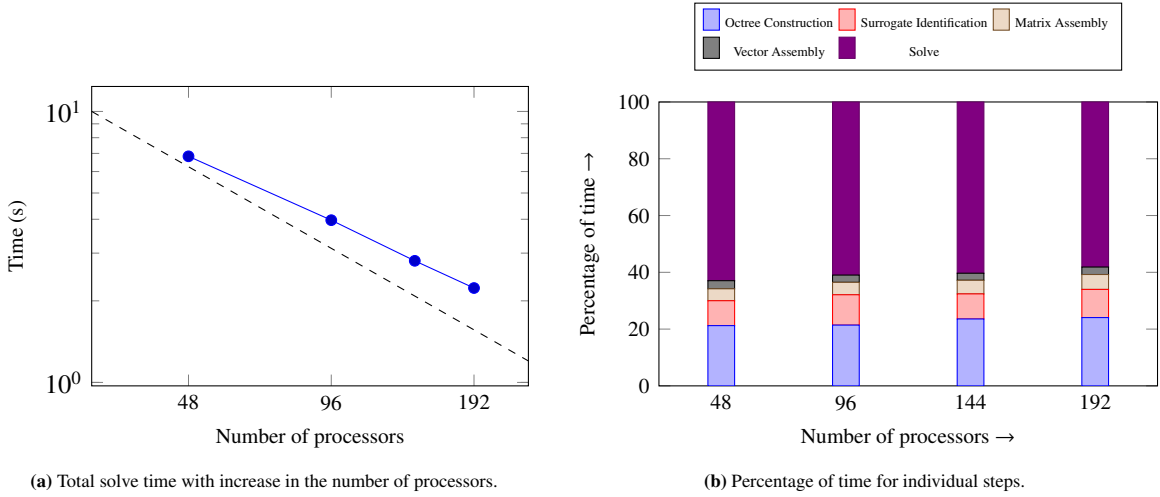
\begin{figure}[t!]
\centering
\begin{subfigure}{.45\linewidth}
    \centering
    \begin{tikzpicture}
    \begin{loglogaxis}[ 
        width=0.99\linewidth,
        height=0.75\linewidth, 
        xtick={48, 96,  192},  
        xticklabel style={rotate=0,font=\footnotesize},
        xticklabels={$48$,$96$, $192$}, 
        scaled y ticks=true,
        xlabel={\footnotesize Number of processors},
        ylabel={\footnotesize Time (s)},
        legend style={at={(0.37,0.2)}, anchor=north, nodes={scale=0.65, transform shape}}, 
        legend columns=3,
        xmin=30,
        xmax=250	
        ]
        \addplot table [x={proc},y={T_0p5_total},col sep=comma] {Figures/Scaling/ScalingCircle.txt};
        \addplot +[mark=none, black, dashed] [domain=30:250]{300/x};
    \end{loglogaxis}
    \end{tikzpicture}
    \caption{Total solve time with increase in the number of processors.}
    \label{fig:ScalingCircle_Total}
\end{subfigure}
\begin{subfigure}{.5\linewidth}
    \centering
    \begin{tikzpicture}
    \begin{axis}[
        width=0.99\linewidth,
        height=0.65\linewidth,
        ybar stacked,
        ymin=0,
        ymax=100,
        xticklabel style={rotate=0,font=\footnotesize},
        xticklabels = {$48$,$96$,$144$,$192$}, %
        xtick=data, 
        enlarge x limits={abs=0.5},
        ylabel={\footnotesize Percentage of time $\rightarrow$},
        xlabel={\footnotesize Number of processors $\rightarrow$},
        yticklabel style={rotate=0,font=\footnotesize},
        legend style={at={(0.5,1.35)},anchor= north,legend columns=3}, 
        ]
        \addplot table[x expr=\coordindex,y expr=\thisrow{T_0p5_OctreeCreate}*100/\thisrow{T_0p5_total},col sep=comma] {Figures/Scaling/ScalingCircle.txt};
        \addplot table[x expr=\coordindex,y expr=\thisrow{T_0p5_SurrogateIndentify}*100/\thisrow{T_0p5_total},col sep=comma] {Figures/Scaling/ScalingCircle.txt};
        \addplot table[x expr=\coordindex,y expr=\thisrow{T_0p5_MatAssembly}*100/\thisrow{T_0p5_total},col sep=comma] {Figures/Scaling/ScalingCircle.txt};
        \addplot table[x expr=\coordindex,y expr=\thisrow{T_0p5_VecAssembly}*100/\thisrow{T_0p5_total},col sep=comma] {Figures/Scaling/ScalingCircle.txt};
        \addplot table[x expr=\coordindex,y expr=(\thisrow{T_0p5_total}-\thisrow{T_0p5_VecAssembly} - \thisrow{T_0p5_MatAssembly} - \thisrow{T_0p5_SurrogateIndentify} - \thisrow{T_0p5_OctreeCreate})*100/
        \thisrow{T_0p5_total},col sep=comma] {Figures/Scaling/ScalingCircle.txt};
        \legend{\tiny Octree Construction, \tiny Surrogate Identification, \tiny Matrix Assembly, \tiny Vector Assembly, \tiny Solve}
    \end{axis}
    \end{tikzpicture}
    \caption{Percentage of time for individual steps.}
    \label{fig:ScalingCircle_percentage}
\end{subfigure}   
\caption{Scaling behavior and percentage of time at different stages of the SBM computation in 2D on TACC ~\Stampede.} 
\label{fig:ScalingCircle}
\end{figure}

\subsection{Parallel computing scaling}
\label{subsec:ComplexGeo}

\begin{figure}[t!]
\centering
\begin{subfigure}{.45\linewidth}
    \centering
    \begin{tikzpicture}
    \begin{loglogaxis}[ 
        width=0.99\linewidth,
        height=0.75\linewidth, 
        xtick={320,640,1280},  
        xticklabel style={rotate=0,font=\footnotesize},
        xticklabels={$320$,$640$, $1280$}, 
        scaled y ticks=true,
        xlabel={\footnotesize Number of processors},
        ylabel={\footnotesize Time (s)},
        legend style={at={(0.37,0.2)}, anchor=north, nodes={scale=0.65, transform shape}}, 
        legend columns=3,
        xmin=300,
        xmax=1800,
        ymin=30,
        ymax=400,
        ]
        \addplot table [x={proc},y={T_0p5_total},col sep=comma] {Figures/Scaling/ScalingBunny.txt};
        \addplot +[mark=none, black, dashed] [domain=300:1800]{70000/x};
    \end{loglogaxis}
    \end{tikzpicture}
    \caption{Total solve time with increase in the number of processors.}
    \label{fig:ScalingBunny_Total}
\end{subfigure}
\begin{subfigure}{.5\linewidth}
    \centering
    \begin{tikzpicture}
    \begin{axis}[
        width=0.99\linewidth,
        height=0.65\linewidth,
        ybar stacked,
        ymin=0,
        ymax=100,
        xticklabel style={rotate=0,font=\footnotesize},
        xticklabels = {$320$,$640$,$960$,$1280$,$1600$}, %
        xtick=data, 
        enlarge x limits={abs=0.5},
        yticklabel style={rotate=0,font=\footnotesize},
        ylabel={\footnotesize Percentage of time $\rightarrow$},
        xlabel={\footnotesize Number of processors $\rightarrow$},
        legend style={at={(0.5,1.35)},anchor= north,legend columns=3}, 
        ]
        \addplot table[x expr=\coordindex,y expr=\thisrow{T_0p5_OctreeCreate}*100/\thisrow{T_0p5_total},col sep=comma] {Figures/Scaling/ScalingBunny.txt};
        \addplot table[x expr=\coordindex,y expr=\thisrow{T_0p5_SurrogateIndentify}*100/\thisrow{T_0p5_total},col sep=comma] {Figures/Scaling/ScalingBunny.txt};
        \addplot table[x expr=\coordindex,y expr=\thisrow{T_0p5_MatAssembly}*100/\thisrow{T_0p5_total},col sep=comma] {Figures/Scaling/ScalingBunny.txt};
        \addplot table[x expr=\coordindex,y expr=\thisrow{T_0p5_VecAssembly}*100/\thisrow{T_0p5_total},col sep=comma] {Figures/Scaling/ScalingBunny.txt};
        \addplot table[x expr=\coordindex,y expr=(\thisrow{T_0p5_total}-\thisrow{T_0p5_VecAssembly} - \thisrow{T_0p5_MatAssembly} - \thisrow{T_0p5_SurrogateIndentify} - \thisrow{T_0p5_OctreeCreate})*100/
        \thisrow{T_0p5_total},col sep=comma] {Figures/Scaling/ScalingBunny.txt};
        \legend{\tiny Octree Construction, \tiny Surrogate Identification, \tiny Matrix Assembly, \tiny Vector Assembly, \tiny Solve}
    \end{axis}
    \end{tikzpicture}
    \caption{Percentage of time for individual steps.}
    \label{fig:ScalingBunny_percentage}
\end{subfigure}   
\caption{Scaling behavior and percentage of time at different stages of the SBM computation for 3D Stanford Bunny on TACC ~\Stampede.} 
\label{fig:ScalingBunny}
\end{figure}
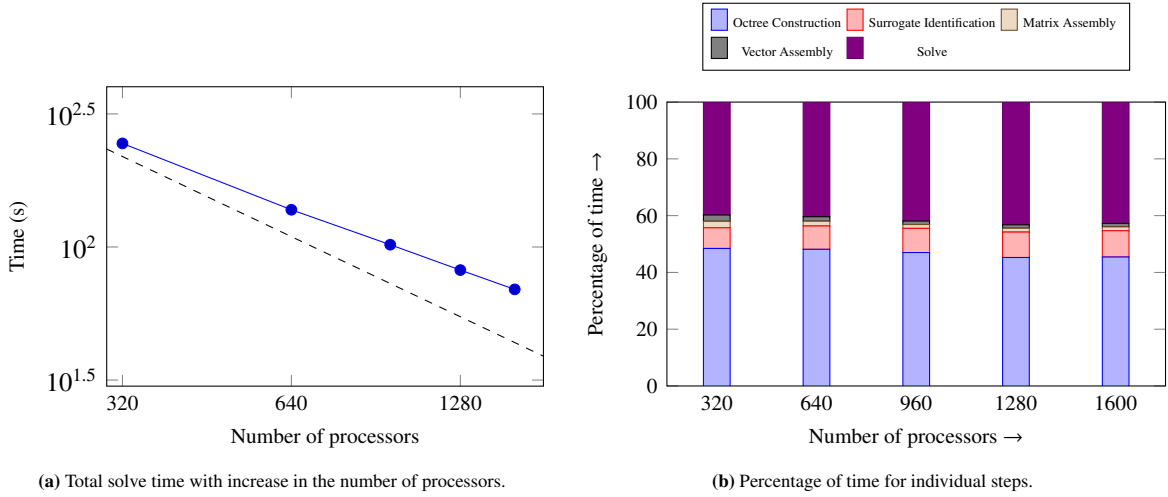

Finally, we present some scaling results of our framework on  TACC~\Stampede ~SKX and ~ICX nodes. To perform the scaling test, we consider the problem described in \secref{subsec:PoissonDisk}. We choose  $\lambda$ of 0.5 for the scaling test as it conforms to the \emph{optimal} surrogate. The strong scaling test considers octree at a level of 10 ($h = 1/2^{10}$). \figref{fig:ScalingCircle_Total} shows the variation of total solve time with increasing the number of processors.; with a near-ideal scaling behavior. \figref{fig:ScalingCircle_percentage} shows the percentage of time taken by the different stages. We observe that the extra step for constructing the surrogate amounts to almost 10\% of the total time. We note that this step, along with octree construction, needs to be performed only once for a static mesh and will be amortized to a smaller fraction for transient problems. As expected, we observe that the overall runtime is dominated by linear algebra solve.

In addition to the 2D case, we conducted a scaling study on a three-dimensional Stanford Bunny, as presented in \secref{subsec:ComplexGeo}. For the strong scaling test, we utilized an octree level of 9 ($h = 1/2^{9}$). The total solving time with respect to the number of processors is shown in \figref{fig:ScalingBunny_Total}, while \figref{fig:ScalingBunny_percentage} illustrates the percentage of time taken by each stage. Unlike the 2D case, the incomplete octree construction time became the problem's bottleneck. However, the octree construction procedure is a one-time event if we encounter time-dependent problems. Notably, the matrix assembly time is similar to the vector assembly time for this three-dimensional case, given that we calculate and store distance function information during vector assembly and reuse it during matrix assembly to minimize the distance function calculation time. Furthermore, we implemented k-D tree to improve the time required to calculate the distance function in three-dimensional complex shapes. Thus, our algorithms and implementation ensure good scalability, making this approach a practical strategy for solving PDEs in complex domains using conceptually simple octree meshes.

\section{Conclusions and Future Work}
\label{Sec:Conclusions}

By shifting the enforcement of boundary conditions from the actual boundary to a surrogate boundary, the SBM allows for the use of Cartesian meshes, eliminating the need for laborious and time-consuming body-fitted meshes around complex geometries. In this work, we answer some key questions regarding the scalable and accurate deployment of the Shifted Boundary Method. The key findings of this work are as follows: (a) identification of an optimal surrogate boundary parameter that greatly reduces numerical error in the SBM, (b) rigorous theoretical analysis demonstrating the optimal convergence of SBM on extended surrogate domains, (c) successful deployment of the SBM on massively parallel octree meshes, including handling of incomplete octrees, and (d) successful application of the SBM to various simulations involving complex shapes, including those with sharp corners and different topologies, with a focus on Poisson's equation and linear elasticity equations. This work sets the stage for a massively parallel, octree-based, general-purpose solution framework---using the SBM---for solving PDEs on arbitrarily complex geometries. 

There are several avenues for future developments. One avenue we are actively exploring involves extending the SBM on the octree framework to multi-physics and coupled PDEs, including Navier-Stokes, Cahn-Hilliard Navier-Stokes, and the Possion-Nernst-Planck equations. Another avenue is to extend the framework to account for moving boundaries, with a natural extension to efficiently model fluid-structure interaction problems across complex geometries. Another active avenue of research is to develop robust preconditioners and architecture-aware solvers (for example, GPU-accelerated multigrid methods) for such SBM on octree approaches. A final straightforward extension is to explore the utility of higher-order basis functions and their tradeoff on error vs. time-to-solve.

\section*{Acknowledgements}
This work was partly supported by the National Science Foundation under the grants NSF LEAP-HI 2053760, NSF CNS 1954556, and NSF OAC 1750865. BG, AK, KS, and CHY are supported in part by AI Research Institutes program supported by NSF and USDA-NIFA under AI Institute: for Resilient Agriculture, Award No. 2021-67021-35329. Guglielmo Scovazzi has been supported by National Science Foundation under Grant 2207164, Division of Mathematical Sciences (DMS).

\bibliographystyle{elsarticle-num-names.bst}
\bibliography{SBMBibs}
\clearpage
\appendix
\section{Proof of Proposition~\ref{prop:extension}}
\label{app:proof_ext}

 \noindent \textbf{Proposition.}
 There exists an extension $\bar{u}$ of $u$ in $\Oe$, such that:   
\begin{enumerate}
    \item[a)] $-\Delta \bar{u} = \bar{f}$ in $\Oe$; and 
    \item[b)] if $u \in H^2(\Omega)$, then $\bar{u} \in H^2(\tO)$, with 
    $\| \bar{u} \|_{H^2(\tO)} \leq C \; \|u\|_{H^2(\Omega)}$. 
\end{enumerate}

\noindent {\it Proof}

To fulfill a) and b) above, let $\Omega_{\Delta} = \Oe \setminus \mathrm{cl}(\Omega)$, where $\mathrm{cl}(\Omega)$ is the closure of $\Omega$.
The boundary $\partial \Omega_{\Delta}$ of the set $\Omega_{\Delta}$ is decomposed as the union of two disjoint boundaries, namely $\partial \Omega_{\Delta} = \Gamma_0 \cup \Gamma_1$, with $\Gamma_0=\partial \Omega$ and $\Gamma_1 = \partial \Oe$.
Then we define
\begin{equation}
\label{eq:ext_u}
    \bar{u} \; = \;
    \begin{cases}
    u \quad \mbox{in } \mathrm{cl}(\Omega) \; , \\
    \hat{u} \quad \mbox{in } \Omega_{\Delta} \; , 
    \end{cases}
\end{equation}
where $\hat{u}$ is the solution of 
\begin{equation}
\label{eq:ext_u1}
   \begin{cases}
    -\Delta \hat{u} = \bar{f} \quad \mbox{in } \Omega_{\Delta} \; ,\\
    \quad \, \, \hat{u} = u_D \quad \mbox{on } \Gamma_0 \; , \\
    \quad \, \, \hat{u} = \varphi \quad \mbox{on } \Gamma_1 \; ,
    \end{cases}
\end{equation}
and $\varphi$ is chosen so that
\begin{equation}
\label{eq:norm_match}
    \frac{\partial \hat{u}}{\partial n}
    = 
    \frac{\partial u}{\partial n}
    \quad \mbox{on } \Gamma_0 \; , \\  
\end{equation}
with $\partial (\cdot) / \partial n$ the normal derivative along $\boldsymbol{n}$ (i.e.\, along the unit normal to $\Gamma_0$ pointing outside of $\Omega$). 
To motivate~\eqnref{eq:norm_match}, observe that $\Gamma_0$ and $\Gamma_1$ are smooth curves and that $u_D \in H^{3/2}(\Gamma_0)$.
Since $f \in L^2(\Omega)$, then, by the regularity theorem of elliptic problems, we have $u \in H^2(\Omega)$ and consequently $\displaystyle \frac{\partial u}{\partial n} \in H^{1/2}(\Gamma_0)$. 
Similarly, if $\varphi \in H^{3/2}(\Gamma_1)$, then from $\bar{f} \in L^2(\Omega_{\Delta})$ and $u_D \in H^{3/2}(\Gamma_0)$ we deduce $\hat{u} \in H^2(\Omega_{\Delta})$ and consequently $\displaystyle \frac{\partial \hat{u}}{\partial n} \in H^{1/2}(\Gamma_0)$.
In order to deduce $\bar{u} \in H^2(\Oe)$, we need to have 
$$
u_{|\Gamma_0} = \hat{u}_{|\Gamma_0} \; ,
$$
which is true since $u_{|\Gamma_0} = \hat{u}_{|\Gamma_0} = u_D$ on $\Gamma_0$, and
$$
\nabla u_{|\Gamma_0} = \nabla \hat{u}_{|\Gamma_0} \; .
$$
This last condition can be decomposed into the matching of the component of the gradient normal to $\Gamma_0$, which is precisely~\eqnref{eq:norm_match}, and the matching of the component of the gradient tangent to $\Gamma_0$, namely
\begin{equation}
    \frac{\partial \hat{u}}{\partial \tau}
    =
    \frac{\partial \hat{u}}{\partial \tau}
    \quad \mbox{on } \Gamma_0 \; , \\  
\end{equation}
which is true since they both coincide with $\displaystyle \frac{\partial u_D}{\partial \tau}$ on $\Gamma_0$.

Next, we verify the existence of $\varphi$ in~\eqnref{eq:ext_u1} such that~\eqnref{eq:norm_match} is verified.  
The solution of~\eqnref{eq:ext_u1} is equivalent to the solution of the two sub-problems 
\begin{subequations}
\begin{align}
\label{eq:ext_u1_sub1}
   &\begin{cases}
    -\Delta \hat{u}_1 = \bar{f} \quad \ \ \mbox{in } \Omega_{\Delta} \; ,\\
    \quad \, \, \hat{u}_1 = u_D \quad \mbox{on } \Gamma_0 \; , \\
    \quad \, \, \hat{u}_1 = 0 \quad \ \ \mbox{on } \Gamma_1 \; ,
    \end{cases}
    \; \, \\[.1cm]
\label{eq:ext_u1_sub0}
   &\begin{cases}
    -\Delta \hat{u}_0 = 0 \quad \mbox{in } \Omega_{\Delta} \; ,\\
    \quad \, \, \hat{u}_0 = 0 \quad \mbox{on } \Gamma_0 \; , \\
    \quad \, \, \hat{u}_0 = \varphi \quad \mbox{on } \Gamma_1 \; ,
    \end{cases}
\end{align}
\end{subequations}
when we set $\hat{u}=\hat{u}_1+\hat{u}_0$, so that condition~\eqnref{eq:norm_match} becomes
\begin{equation}
\label{eq:norm_match_sub}
    \frac{\partial \hat{u}_0}{\partial n}
    = 
    \frac{\partial u}{\partial n}
    - \frac{\partial \hat{u}_1}{\partial n}
    \; \in H^{1/2}(\Gamma_0) \; .
\end{equation}
Let us then introduce the the operator 
\begin{align*}
    \mathcal{T} : H^{1/2}(\Gamma_1) &\rightarrow \; H^{-1/2}(\Gamma_0) \\
        \varphi &\mapsto \; \left. \frac{\partial \hat{u}_0}{\partial n}\right|_{\Gamma_0}
\end{align*}
and check that $\mathcal{T}$ is {\it onto} (i.e., {\it surjective}), namely that $\mathrm{Im}(\mathcal{T}) = \left\{ \mathcal{T} \varphi : \varphi \in H^{1/2}(\Gamma_1) \right\}$ coincides with $H^{-1/2}(\Gamma_0)$.
We argue by contradiction: since $\mathrm{Im}(\mathcal{T})$ is a closed subspace of $H^{-1/2}(\Gamma_0)$, assume that its orthogonal space contains elements $\eta \neq 0$. Any such $\eta \in H^{-1/2}(\Gamma_0)$ then satisfies
$$
(\eta , \mathcal{T} \varphi)_{H^{-1/2}(\Gamma_0)} = 0 \; , \qquad \forall \varphi \in H^{-1/2}(\Gamma_0) \; .
$$
Note that the inner product in $H^{-1/2}(\Gamma_0)$ is the duality pairing 
$$
_{H^{-1/2}(\Gamma_0)} \langle \mathcal{T} \varphi \; , \, R \eta \rangle_{H^{1/2}(\Gamma_0)} \; ,
$$
where $R \eta = z_{|\Gamma_1}$ with $z$ the solution of the problem
\begin{align}
\label{eq:z}
   &\begin{cases}
    \; \;-\Delta z = 0 \quad \mbox{in } \Omega_{\Delta} \; ,\\
    \; \; \quad \, \, z = 0 \quad \mbox{on } \Gamma_1 \; , \\[.1cm]
    \; \; \; \; \displaystyle \frac{\partial z}{\partial n} = \eta \quad \mbox{on } \Gamma_0 \; .
    \end{cases}
\end{align}
Hence, for any $\varphi \in H^{1/2}(\Gamma_1)$, we have
\begin{align}
    0
    &= \;
    _{H^{-1/2}(\Gamma_0)} \langle \frac{\partial \hat{u}_0}{\partial n} \; , \, z \rangle_{H^{1/2}(\Gamma_0)} 
    \nonumber \\
    &= \;
    _{H^{-1/2}(\Gamma_0)} \langle \frac{\partial \hat{u}_0}{\partial n} \; , \, z \rangle_{H^{1/2}(\Gamma_0)} 
    \; + \;
    \underbrace{
    _{H^{-1/2}(\Gamma_1)} \langle \frac{\partial \hat{u}_0}{\partial n} \; , \, z \rangle_{H^{1/2}(\Gamma_1)} 
    }_{\mbox{$=0$ by~\eqnref{eq:z}}}  
    \nonumber \\
    &= \;
    _{H^{-1/2}(\partial \Omega_{\Delta})} \langle \frac{\partial \hat{u}_0}{\partial n} \; , \, z \rangle_{H^{1/2}(\partial \Omega_{\Delta})} 
    \nonumber \\
    &= \;
    \int_{\Omega_{\Delta}} \underbrace{\Delta \hat{u}_0}_{\mbox{$=0$ by~\eqnref{eq:ext_u1_sub0}}} \, z
    \; + \;
    \int_{\Omega_{\Delta}} \nabla \hat{u}_0 \cdot \nabla z 
    \qquad \qquad \qquad \qquad \ \ \ \ \mbox{(using integration by parts)}
    \nonumber \\
    &= \;
    \int_{\Omega_{\Delta}} \nabla \hat{u}_0 \cdot \nabla z 
    \nonumber \\
    &= \;
    \int_{\Omega_{\Delta}} \underbrace{\Delta z}_{\mbox{$=0$ by~\eqnref{eq:z}}} \, \hat{u}_0
    \; + \;
    _{H^{-1/2}(\partial \Omega_{\Delta})} \langle \frac{\partial z}{\partial n} \; , \, \hat{u}_0 \rangle_{H^{1/2}(\partial \Omega_{\Delta})}
    \qquad \qquad \mbox{(using integration by parts)}
    \nonumber \\
    &= \;
    \underbrace{
    _{H^{-1/2}(\Gamma_0)} \langle \frac{\partial z}{\partial n} \; , \, \hat{u}_0 \rangle_{H^{1/2}(\Gamma_0)} 
    }_{\mbox{$=0$ by~\eqnref{eq:ext_u1_sub0}}}  
    \; + \;
    _{H^{-1/2}(\Gamma_1)} \langle \frac{\partial z}{\partial n} \; , \, \hat{u}_0 \rangle_{H^{1/2}(\Gamma_1)} 
    \nonumber \\
    &= \;
    _{H^{-1/2}(\Gamma_1)} \langle \eta \; , \, \varphi \rangle_{H^{1/2}(\Gamma_1)} 
    \; .
 \end{align}
 We conclude then that $\eta$ is the null linear form on $H^{1/2}(\Gamma_0)$, that is $\eta=0$.
 As a consequence of $\mathcal{T}$ being {\it onto}, picking
 $$
 \eta
 \;=\;
\frac{\partial u}{\partial n}
-
\frac{\partial \hat{u}_1}{\partial n}
\; ,
$$
there exists $\phi \in H^{1/2}(\Gamma_1)$ such that the solution of~\eqnref{eq:ext_u1_sub0} satisfies
$$
\frac{\partial \hat{u}_0}{\partial n} = \eta \; .
$$
In addition, since $\eta \in H^{1/2}(\Gamma_0)$, we get $\hat{u}_0 \in H^2(\Omega_{\Delta})$ and $\varphi = (\hat{u}_0)_{| \Gamma_1} \in H^{3/2}(\Gamma_1)$, whence $\hat{u} \in H^2(\Omega_\Delta)$ as desired. $\Box$

\end{document}